\documentclass[letterpaper,12pt,peerreviewca,draftcls]{IEEEtran}
\usepackage{csm16}
\usepackage[margin=1in]{geometry}
\usepackage[nolists,nomarkers,tablesfirst]{endfloat} 
\usepackage{amsmath,amsfonts,amssymb,mathtools} 
\usepackage{tikz}

\usepackage[labelformat=simple]{subcaption}

\usepackage{algorithmic}
\usepackage{url}
\usepackage{graphicx,xcolor}
\usepackage{verbatim}
\makeatletter
\newcommand{\verbatimfont}[1]{\def\verbatim@font{#1}}%
\makeatother
\verbatimfont{\ttfamily\small}

\newcommand{\bi}{\begin{itemize}}\newcommand{\ei}{\end{itemize}}
\newcommand{\be}{\begin{equation}}\newcommand{\ee}{\end{equation}}
\newcommand{\bee}{\begin{enumerate}}\newcommand{\eee}{\end{enumerate}}
\newcommand{\bea}{\begin{eqnarray}}\newcommand{\eea}{\end{eqnarray}}
\newcommand{\beas}{\begin{eqnarray*}}\newcommand{\eeas}{\end{eqnarray*}}
\newcommand{\bc}{\begin{center}}\newcommand{\ec}{\end{center}}
\usepackage[left,pagewise]{lineno}
\usepackage[english]{babel} 
\usepackage{blindtext}

\def\Re{\mathbb{R}}
\def\R{\mathbb{R}}

\def\argmin{\mathop{\text{\rm arg\,min}}}

\def\notes#1{\marginpar{\tiny #1}\typeout{Notes!
Notes!
Notes!
}}
\renewcommand{\notes}[1]{\typeout{notes!}}
\def\FRAC#1#2#3{\genfrac{}{}{}{#1}{#2}{#3}}
\def\half{{\mathchoice{\FRAC{1}{1}{2}}%
{\FRAC{2}{1}{2}}%
{\FRAC{3}{1}{2}}%
{\FRAC{4}{1}{2}}}}

\newcommand{\tr}{\mbox{tr}}

\def\Re{\field{R}}

\def\clL{{\cal L}}
\def\clP{{\cal P}}
\def\clZ{{\cal Z}}

\def\E{{\sf E}}

\def\clC{{\cal C}}

\def\R{\mathbb{R}}

\def\clZ{{\cal Z}}

\newtheorem{theorem}{Theorem}

\newtheorem{definition}{Definition}

\newtheorem{remark}{Remark}
\newtheorem{proposition}{Proposition}

\def\beq{\begin{eqnarray}} 
\def\bc{\begin{center}} 
\def\be{\begin{enumerate}}
\def\bi{\begin{itemize}} 
\def\bs{\begin{small}}
\def\bS{\begin{slide}}
\def\ec{\end{center}} 
\def\ee{\end{enumerate}}
\def\ei{\end{itemize}}
\def\es{\end{small}}
\def\eS{\end{slide}}
\def\eeq{\end{eqnarray}}

\newcommand{\newP}[1]{\medskip\noindent{\bf #1:}}

\newcommand{\ud}{\,\mathrm{d}}

\def\Re{\mathbb{R}}
\def\E{{\sf E}}

\def\argmin{\mathop{\text{\rm arg\,min}}}

\def\clY{{\cal Y}}

\def\Thm#1{Thm.~\ref{#1}}
\def\Prop#1{Prop.~\ref{#1}}

\def\clD{{\cal D}}

\def\clL{{\cal L}}
\def\clP{{\cal P}}
\def\clZ{{\cal Z}}

\renewcommand{\Re}{\mathbb{R}}

\def\FRAC#1#2#3{\genfrac{}{}{}{#1}{#2}{#3}}

\def\clA{{\cal A}}

\def\clC{{\cal C}}

\def\clD{{\cal D}}

\def\clF{{\cal F}}

\def\clL{{\cal L}}

\def\clM{{\cal M}}

\def\clP{{\cal P}}

\def\clU{{\cal U}}
\def\clV{{\cal V}}

\def\clY{{\cal Y}}
\def\clZ{{\cal Z}}
\def\E{{\sf E}}
\def\bS{\mathbb{S}}
\def\sJ{{\sf J}}
\def\bsJ{{\sf J}}
\def\ones{{\sf 1}}
\def\sP{{\sf P}}
\def\tsP{{\tilde{\sf P}}}
\def\tE{{\tilde{\sf E}}}

\def\tp{{\hbox{\rm\tiny T}}}

\def\Nsp{{\sf N}} 
\def\Rsp{{\sf R}} 
\def\hE{\tilde{\sf E}}

\def\sQ{{\sf Q}}

\def\opt{{\text{\rm (opt)}}}

\newcommand{\revblue}[1]{{\color{blue}#1}}

\newcounter{rmnum}

\newenvironment{romannum}{\begin{list}{{\upshape (\roman{rmnum})}}{\usecounter{rmnum}
			\setlength{\leftmargin}{12pt}
			\setlength{\rightmargin}{8pt}
			\setlength{\itemsep}{2pt}
			\setlength{\itemindent}{-1pt}
	}}{\end{list}}

\newcounter{anum}

\title{Arrow of Time in Estimation and Control\\[5pt]
  \Large Duality Theory Beyond the Linear Gaussian Model}

\medskip

\author{Jin Won Kim and Prashant G. Mehta\\[10pt] September 19, 2024}

\newif\ifPDF \ifx\pdfoutput\undefined\PDFfalse \else\ifnum\pdfoutput > 0\PDFtrue \else\PDFfalse \fi \fi
\ifPDF 
\usepackage[pdftex, plainpages = false, colorlinks=true, linkcolor=black, citecolor = green!50!blue, urlcolor = blue, filecolor=black, pagebackref=false, hypertexnames=false,  pdfpagelabels ]{hyperref}
\fi

\begin{document}
	\maketitle
	\CSMsetup

The House of Control stands on two pillars: Estimation and Control.
The elegant symmetry of these pillars is evident in Kalman's {\em
  Principle of Duality} (PoD)~\cite[pp.~489]{kalman1960general}:

\noindent Principle of Duality (PoD). Considering the
class of feedback systems with linear control law,
the dual plant [...] is obtained by the
following steps: 
\begin{romannum}
\item Replace [the state transition matrix] by its dual [the transpose],
\item Interchange input and output constraints, 
\item Reverse the direction of time.
\end{romannum}

\begin{figure}
        \centering
        \includegraphics[width=0.45\textwidth]{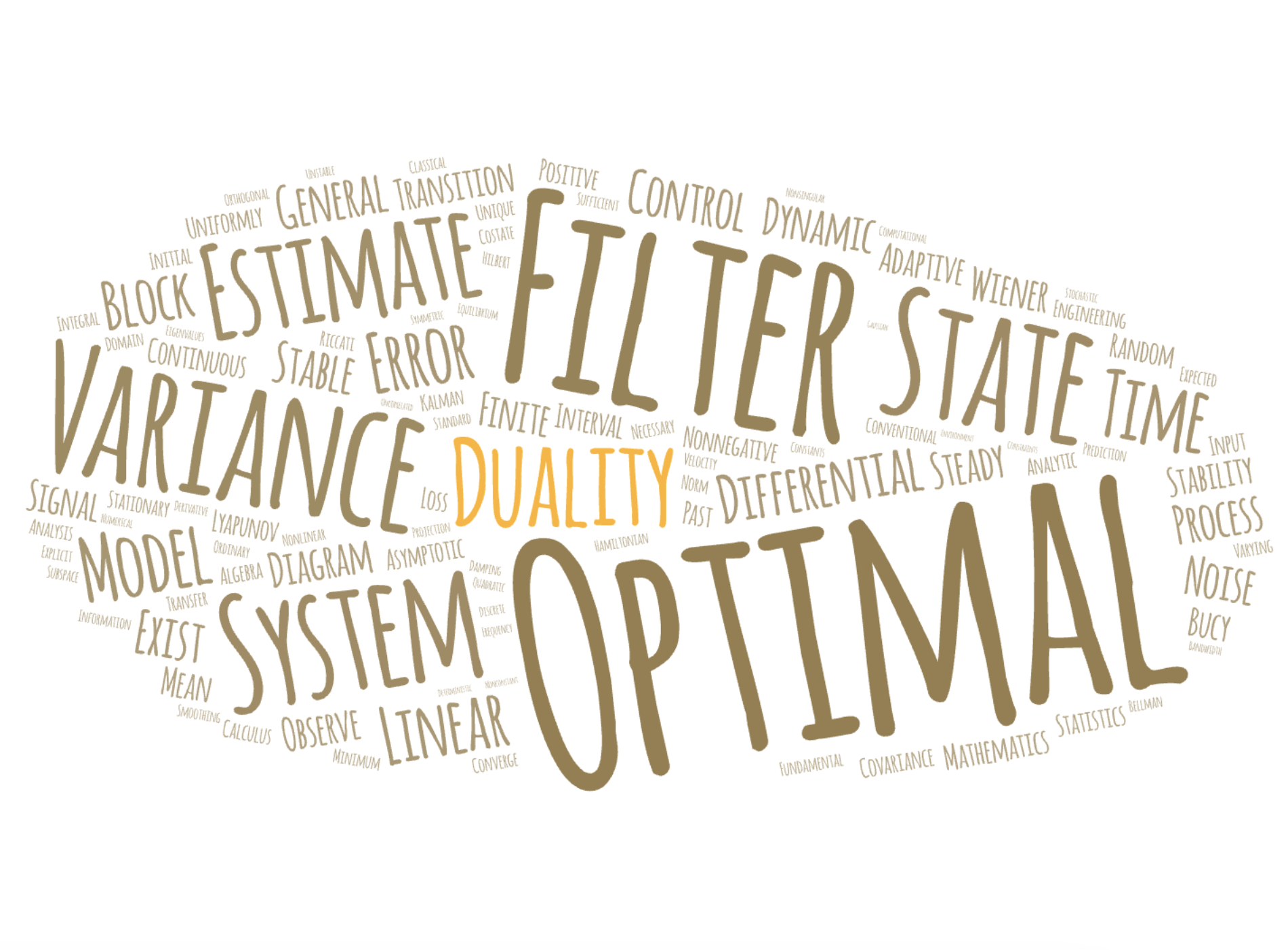}
\caption{Word cloud of the Kalman-Bucy's
  paper~\cite{kalman1961}. The word `dual' or `duality' appears
  a total of 42 times in the paper.  The word `optimal' with 78
  appearances is the highest frequency word in the paper.  The commonly used English language words
such as 'and', 'the' etc., have been excluded from the word cloud.}
\label{fig:word_cloud}
\end{figure}

Duality is coeval with the origin of modern systems and control
theory and it underpins much of
Kalman's early
work.  For example, his celebrated paper together with Bucy describing
the Kalman-Bucy filter draws extensively on the PoD -- between the optimal filter and a certain linear quadratic (LQ) optimal control problem (see word cloud in Fig.~\ref{fig:word_cloud})~\cite{kalman1961}.  Notably,
duality help explains why the
Riccati equation is the fundamental equation for \textit{both} optimal
filtering and optimal control.

Duality is
expressed in two inter-related manners:
\begin{romannum}
\item Duality between observability and controllability; and
\item Duality between optimal filtering and optimal control.
\end{romannum}
The second of these items means expressing one type of problem as
another type of problem.  In the House of Control, the main interest
is to solve an optimal filtering (or estimation) problem by expressing
it as an optimal control problem.

Sixty-five years have elapsed since the PoD was first stated.  Duality
is now the foundation upon which the two pillars of the House of
Control stand. As part of most control curricula, controllability
and observability are introduced as dual concepts in the first
graduate course.  Duality permeates every
facet of the subject's instruction: (i) \textbf{Analysis}: dual nature of rank
tests, Gramians, Lyapunov functions, and their relevance to minimal
realization and to Hankel model reduction; and (ii) \textbf{Design}:
synthesis of stabilizing controller and Luenberger observer are dual
problems.

Our objective in this paper is to describe extensions of duality
theory for nonlinear stochastic systems (hidden Markov models or
HMMs).  We argue that, despite six decades of work on this topic,
duality theory is incomplete.  We discuss both the difficulty in
extending duality to HMMs, as well as its recent resolution by our
research group at Illinois.



As early as 1977, Hermann
and Krener wrote a pioneering paper ``Nonlinear
controllability and observability'' where it is noted that 
``{\em [nonlinear] duality between “controllability” and “observability”
  [...] is, mathematically, just the duality between vector fields and
  differential forms}''~\cite{hermann1977nonlinear}.  In recent decades, detectability of a
nonlinear system is understood as an ``output-to-state stability (OSS)''
property.  Its definition was originally motivated by Sontag and Wang as follows:  
``{\em Given the central role often played in control theory by the
  duality between input/state and state/output behavior, one may
  reasonably ask what concept obtains if outputs are used instead of
  inputs in the [input-to-state stability (ISS)] definition. [..] it would appear that this dual property, which
  we will call output-to-state stability (OSS), is a natural candidate
  as a concept of nonlinear (zero-) detectability}''\cite{sontag1997output}. A variant of the
OSS, the so called incremental output-to-state stability (i-OSS), is 
a standard assumption to establish stability of a class of nonlinear 
estimators known as the moving horizon estimators
(MHE).
A historical account appears in \cite[Sec. 4.7]{rawlings2017model}, where
Rawlings et al. note ``{\em establishing duality
[of the optimal estimator] with the optimal regulator is a
favorite technique for establishing estimator stability}''.


Given all of this history, how can it be that duality is incomplete or
not understood for HMMs?  The issue is not that the HMMs are
stochastic while the OSS definitions and MHE algorithms are for
deterministic systems. 
Even though MHE is often introduced and studied in deterministic settings, the algorithm has its roots in the study of
stochastic linear Gaussian systems.  Specifically, MHE is based on the so called
minimum energy duality which originated in the early work (in 1960s) of Bryson
and Frazier, Mortensen, and the important PhD thesis of Hijab
(the 1980 PhD thesis was supervised by
Krener)~\cite{bryson1963smoothing,mortensen1968,hijab1980minimum}. 
Soon after Hijab's PhD thesis, 
Fleming and Mitter wrote an influential paper where they 
describe the log
transformation link between the Zakai equation of nonlinear
(stochastic) filtering and the Hamilton-Jacobi-Bellman (HJB) equation of nonlinear stochastic optimal control~\cite{fleming1982optimal}.  
One consequence is that the negative log of the conditional probability
(posterior) 
density for the filtering problem is a
value function for some stochastic optimal control problem. 
Even though the log transform link was noted, the dual optimal control problem itself
was not clarified in~\cite{fleming1982optimal} or the studies that
followed~\cite{fleming1997deterministic,benevs1983relation,pavon1989stochastic,james1988nonlinear}.  The dual optimal control problem was presented for the first time in a 2003 paper from Mitter and Newton~\cite{mitter2003}.
Their construction is based on a variational formulation of the Bayes'
formula in terms of the Kullback-Leibler (K-L) 
divergence. 
The KL-based variational form of Bayes is actually quite classical (see~\cite[Sec.~3]{mitter2000duality},~\cite[Sec.~2.2]{van2006filtering}) and
Mitter-Newton's work may be regarded as a generalization of
the minimum energy duality of Mortensen-Hijab. 
A notable ensuing contribution is the PhD thesis of van
Handel where the Fleming-Mitter-Newton duality is used to
obtain results on nonlinear filter stability~\cite{van2006filtering}. 

In this paper, it is argued that the duality
between estimation and control is incomplete in part because of the
issue of time reversal -- step (iii) in Kalman's
PoD.  That there is a problem in generalizing the classical 
duality of Kalman to nonlinear settings has previously been noted in literature.  For
example, in his 2008
paper~\cite{todorov2008general}, Todorov writes: ``{\em Kalman's duality has been known for half a century and has attracted a lot of attention. If a straightforward
generalization to non-LQG settings was possible it would have been
discovered long ago. Indeed we will now show that Kalman's duality,
although mathematically sound, is an artifact of the LQG setting.}''  Todorov attributes the problem to the lack of linearity,
e.g., linear algebraic operations such as the use of matrix transpose
are not applicable to nonlinear systems.  In this paper, 
it is shown that the issue is not so much the lack of linearity but
time reversal. Some of our own recent work on this topic is described where we have presented a solution to this important gap in literature.  

Several comparisons are provided
to illustrate the
differences with prior work as well as
discuss the significance of these differences.  Specifically, it
is shown that the elementary duality between observability and
controllability is an example of Kalman's duality that our work
helps generalize.  Kalman's duality is shown to be distinct from the
minimum energy duality which is implicit in much of the other work
including MHE, papers of Mortensen and Hijab, the log transform of
Mitter-Fleming, and the dual optimal control problem of Mitter-Newton.

In the present paper,
Kalman's duality is referred to as the {\em minimum variance duality}.  The
two terms ``minimum variance duality'' and ``minimum energy duality''
are borrowed from Bensoussan's writings on the subject (which have inspired our own work including the present article).  An 
impetus comes from a curious dichotomy in his
1992 textbook~\cite{bensoussan1992stochastic} on partially observed
systems: While optimization and optimal control techniques are prominently
used for the linear Gaussian model, these are conspicuous by their
absence in the derivation and analysis of the nonlinear filter and
smoothing equations (see~\cite[Sec.~4.8.]{bensoussan1992stochastic}).

\newpage



\section{Arrow of time in estimation and control}

While duality is introduced as part of the first graduate class in
most control curricula, the ``step (iii)'' {\em reversing the direction of
time} in Kalman's PoD is rarely stressed.  In this section,
time reversal is 
illustrated with a discussion of two elementary concepts related to
deterministic and stochastic (linear Gaussian) systems. In the
accompanying sidebar ``Arrow of time in ordinary and stochastic differential equations'', elementary theory of forward-in-time and
backward-in-time differential equations is presented
in a self-contained manner.  This theory explains why the problem of
reversing time is not especially pertinent in deterministic settings but
non-trivial in stochastic settings.     

\subsection{Time reversal for elementary duality between
  observability and controllability}
\label{sec:cntrb_obsvb}

Consider the linear state-output system:
\begin{subequations}\label{eq:LTI-obs}
	\begin{align}
		\frac{\ud x_t}{\ud t} &= A^\tp x_t,\quad x_0 = \xi \label{eq:LTI-obs-a}\\
		z_t &= H^\tp x_t \label{eq:LTI-obs-b}
	\end{align}
\end{subequations}
Over a fixed time interval $[0,T]$, the output is denoted by $z =
\{z_t\in\Re^m:0\le t\le T\}$.  The output $z$ is an element of the
function space $L^2([0,T];\Re^m)$.  

The output depends upon the initial condition $\xi$. This
dependence is indicated by using the superscript whereby the output
$z$ with initial condition $x_0=\xi$ is denoted by $z^{\xi}$.  A basic estimation problem is to determine the initial condition $\xi$ from the output 
$z^{\xi}$. Regarding the solution of this problem, the following definition 
naturally arises.

\begin{definition}\label{def:LTI-obs-def}
	The linear system~\eqref{eq:LTI-obs} is \emph{observable} if:
	\[
	z^{\xi_1} = z^{\xi_2} \quad\Longrightarrow\quad \xi_1 = \xi_2,\quad \forall\,\xi_1,\xi_2\in\Re^d
	\]
\end{definition}

Consider next the following linear input-state system:
\begin{subequations}\label{eq:LTI-ctrl}
\begin{align}
	-\frac{\ud y_t}{\ud t} &= A y_t + H u_t , \quad 0\leq t\leq
                                 T\\
y_T&=0
\end{align}
\end{subequations}
where $u = \{u_t\in\Re^m:0\le t \le T\}$ is the control
input.  Because the terminal condition at final time $t=T$ is
specified, the system~\eqref{eq:LTI-ctrl} is an example of a backward-in-time
control system whose solution at time $t$ is obtained by integrating
backward-in-time from $T$ to $t$ as follows:
\[
y_t = \int_t^T e^{A(s-t) } H u_s \ud s,\quad 0\leq t\leq T
\]
The control $u\in L^2([0,T];\Re^m)$.  

Viewed as a control system, a basic control problem is to design an input $u$ that steers the system from 
a given initial condition $y_0 = \eta$ to $y_T = 0$. Regarding the solution of 
this problem, the following definition naturally arises.

\begin{definition}
	For the linear system~\eqref{eq:LTI-ctrl}, the \emph{controllable subspace} is defined by:
	\[
	\clC = \big\{ \eta \in \Re^d: \exists u \in L^2([0,T];\Re^m) \text{ such that the solution to~\eqref{eq:LTI-ctrl} satisfies } y_0 = \eta \text{ and } y_T= 0\big\}
	\]
	The linear system~\eqref{eq:LTI-ctrl} is \emph{controllable} if $\clC = \Re^d$
\end{definition}

Now, it is well known that observability of the state-output
system~\eqref{eq:LTI-obs} and the controllability of the input-state
system~\eqref{eq:LTI-ctrl} are dual properties.  This nature of this duality
is a special case of the abstract duality in linear algebra
which is reviewed in the sidebar ``Mathematical foundations of
duality''.  The explicit calculations for~\eqref{eq:LTI-obs}
and~\eqref{eq:LTI-ctrl} are described next. 
 
The first step is to define the function spaces which are of the
following two types:
\begin{romannum}
\item \textbf{For the state (initial condition)}. Euclidean space
  $\Re^d$ equipped with the inner-product 
		\[
		\langle \eta,\xi\rangle_{\Re^d} = \eta^\tp \xi,\quad \eta,\xi \in \Re^d
		\]
\item \textbf{For the signal (input and output)}.  The Hilbert space
  of square integrable $\Re^m$-valued signals
  $L^2\big([0,T];\Re^m\big)$ equipped with the inner product
		\[
		\langle u,z\rangle_{L^2([0,T];\Re^m)} = \int_0^T u_t^\tp z_t \ud t,\quad u,z\in L^2\big([0,T];\Re^m\big)
		\]	
\end{romannum}

\begin{figure}
	\centering
	\includegraphics[width=0.85\textwidth]{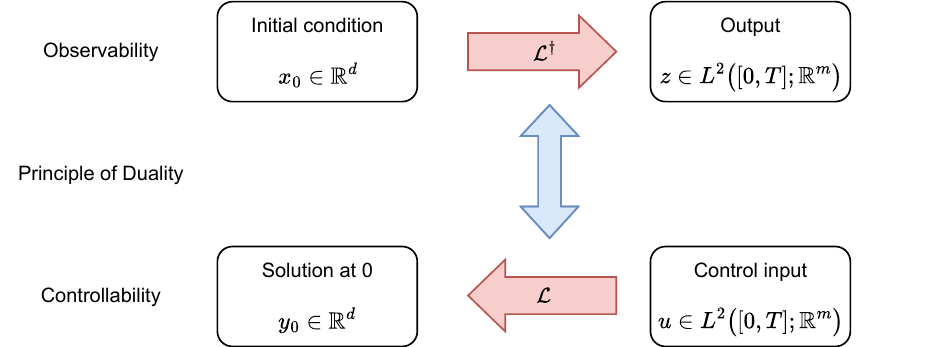}
			\caption{Kalman's duality between controllability and observability: The two are properties of a linear operator and its adjoint.}
		\label{fig:duality_LTI}
	\end{figure}

\newP{Duality} 
For the state-output system~\eqref{eq:LTI-ctrl}, the solution map $u\mapsto y_0$ is used to define a linear operator $\clL: L^2([0,T];\Re^m)\to \Re^d$ as follows:
\begin{equation*}\label{eq:LTI-adjoint}
	\clL u := y_0 = \int_0^T e^{At}Hu_t \ud t
\end{equation*}
The controllable subspace $\clC=\Rsp(\clL)$, the range space of
the operator $\clL$. The operator's adjoint is given by
\[
(\clL^\dagger \xi)(t) = H^\tp e^{A^\tp t}\xi,\quad 0\le t \le T
\]
and represents the solution map from initial condition $\xi\mapsto
z^\xi$ for the state-output system~\eqref{eq:LTI-obs}.  The null-space
of $\clL^\dagger$ is precisely the subspace of initial conditions that
produce zero output, and as such represent the ``unobservable''
directions. 
The duality relationship is given by
\begin{equation*}\label{eq:dual_LTI}
\langle \xi,\clL u\rangle_{\Re^d} = \langle \clL^\dagger \xi, u \rangle_{L^2([0,T];\Re^m)},\quad \forall\, \xi \in \Re^d,\; u\in L^2([0,T];\Re^m)
\end{equation*}
Because of this relationship, the input-state
system~\eqref{eq:LTI-ctrl} is referred to as the \emph{dual control system} for the
state-output system~\eqref{eq:LTI-obs}.  An important implication of
duality, obtained as a corollary of the closed range theorem,  
is that~\eqref{eq:LTI-obs} is observable iff the dual control
system~\eqref{eq:LTI-ctrl} is controllable. 
Figure~\ref{fig:duality_LTI} depicts this duality relationship.

We close by noting the application of Kalman's PoD:
\begin{romannum}
\item The matrices $(A^\tp,H^\tp)$ in
  the state-output system~\eqref{eq:LTI-obs} are replaced by
  their respective transposes $(A,H)$ in the input-state system~\eqref{eq:LTI-ctrl},
\item The input and the output constraints are interchanged.  The
  constraint here is 
  that these are both elements of the common function space $L^2([0,T];\Re^m)$.
\item The direction of time is reversed.  While the state-output
  system~\eqref{eq:LTI-obs} is integrated forward-in-time, the
  input-state system~\eqref{eq:LTI-ctrl} is integrated backward-in-time.
\end{romannum}

\revblue{
\begin{remark}
  Without consideration of inputs and outputs ($H=0$), the system~\eqref{eq:LTI-ctrl} is referred to as an {\em adjoint equation} to system~\eqref{eq:LTI-obs}~\cite{kouba2020adjoint}. Concerning the adjoint equation, Brockett writes~\cite[pp.~44]{brockett2015finite}: ``{\em [..] the key property of adjoint equation is that it propagates the solution of the original equation backwards in time}.'' Therefore, the steps (i) and (iii) of Kalman's PoD may be regarded as the same step based on construction of the adjoint equation. 
\end{remark}

\begin{remark}
  All modern textbooks on linear systems theory include a section or more on duality between controllability and observability.  For works that specifically focus on estimation and include a discussion of duality, see~\cite[Appendix A.17]{rawlings2017model},~\cite[Appendix B.2]{rao_thesis},~\cite[Sec.~15.3.4]{kailath2000linear},~\cite[Ch.~9]{xiong2008introduction},~\cite[Ch.~8]{junkins}. Both~\cite{rawlings2017model,rao_thesis} cite the classic textbook~\cite[Ch.~8]{callier2012linear} for background on duality.
\end{remark}
}

\subsection{Time reversal for the DRE of optimal
  filtering and optimal control}

In stochastic filtering, a basic model is the noisy counterpart of the
state-output system~\eqref{eq:LTI-obs}.  The model is referred to as
the {\em linear Gaussian model} and described by the following It\^{o}
stochastic differential equation (SDE):
\begin{subequations}\label{eq:linear-Gaussian-model}
	\begin{align}
		\ud X_t &= A^\tp X_t \ud t + \sigma \ud B_t,\quad X_0\sim N(m_0,\Sigma_0) \label{eq:linear-Gaussian-model-a}\\
		\ud Z_t &= H^\tp X_t \ud t + \ud W_t,\quad Z_0=0 \label{eq:linear-Gaussian-model-b}
	\end{align}
\end{subequations}
where $X:=\{X_t:t\geq 0\}$
is the $\Re^d$-valued state process, the prior $N(m_0,\Sigma_0)$ is a Gaussian
density with mean $m_0\in \Re^d$ and variance $\Sigma_0 \in
\Re^{d\times d}$,
$Z:=\{Z_t:t\geq 0\}$ is the $\Re^m$-valued observation process, and
$B:=\{B_t:t\geq 0\}$ and
$W:=\{W_t:t\geq 0\}$ are Brownian motions (B.M.) of appropriate dimensions. It is
assumed that $X_0,B,W$ are mutually independent. 
Denote $Q:=\sigma
\sigma^\tp$ and $R$ is the covariance of the observation
noise $W$.  \revblue{It is assumed that $R\succ 0$}.

The filtration
generated by observations up to time $t$ is denoted by
$\clZ_t:=\sigma(\{Z_s:0\leq s\leq t\})$.    
The objective of the linear Gaussian filtering problem is to compute
the conditional expectation 
$\E(f^\tp X_T\mid \clZ_T)$ for a given deterministic $f\in\Re^d$.  

The Kalman-Bucy filter is a recursive algorithm that solves this
problem:
\begin{subequations}\label{eq:KF_eqn}
\begin{align}
\ud m_t  & = A^\tp m_t \ud t + K_t (\ud Z_t - H^\tp m_t \ud t) \label{eq:KF_mean}\\
\frac{\ud \Sigma_t}{\ud t} & = A^\tp \Sigma_t + \Sigma_t A + Q -
                             \Sigma_t H R^{-1} H^\tp \Sigma_t  \label{eq:KF_var}
\end{align}
\end{subequations}
where $K_t = 
\Sigma_t H R^{-1}$ is referred to as the Kalman gain.  
The filter is
initialized with initial condition $(m_0,\Sigma_0)$.  The random
vector $m_t=\E(X_t\mid \clZ_t)$ is the conditional mean, and therefore,
$\E(f^\tp X_T\mid \clZ_T)= f^\tp m_T$.  The deterministic matrix
$\Sigma_t=\E((X_t-m_t)(X_t-m_t)^\tp)$ is the conditional
covariance.  

Related to the PoD, the following remarks appear in the paper of Kalman and Bucy~\cite{kalman1961}:
\begin{romannum}
\item The equation for the covariance $\{\Sigma_t:t\geq 0\}$ is a differential Riccati
  equation (DRE).  DRE is well known as an optimality equation for
  the deterministic LQ optimal control problem. In the standard
  version of the LQ problem, with a forward-in-time linear control system as the
  constraint, the DRE is integrated backward-in-time. 
\item In contrast, the DRE of the Kalman filter is integrated
  forward-in-time.  Therefore, the arrow of time is reversed in going
  from optimal control to optimal filtering, consistent with the
  step~(iii) of Kalman's PoD.
\end{romannum}

The appearance of the DRE suggests that the problem of computing $\E(f^\tp X_T\mid \clZ_T)$ may alternatively be expressed as an optimal control problem.  \revblue{This is indeed the case and the dual optimal control problem is an example of minimum variance duality (the construction appears in the classical textbook of {\AA}str{\"o}m~\cite[Sec.~7.6]{astrom1970})}.  To help motivate the construction, we begin by recalling that the conditional expectation is a solution of an optimization problem:
\[
\E(f^\tp X_T\mid \clZ_T) = \argmin_{S_T\in\clZ_T} \E\big(|f^\tp X_T - S_T|^2\big) 
\]
where ``$S_T\in\clZ_T$'' means that the random variable $S_T$ is $\clZ_T$-measurable (i.e., allowed to depend upon observations made up to time $T$).  The random variable $S_T$ is referred to as an {\em
  estimator} and $\E\big(|f^\tp X_T - S_T|^2\big)$ is the {\em
  mean-squared error (m.s.e)} for the estimator.  The equality means that the conditional expectation $S_T=f^\tp m_T$ is the estimator with the minimum m.s.e..  Because $\Sigma_T$ is the conditional covariance, the minimum m.s.e. is given by
\[
\E\big(|f^\tp X_T - f^\tp m_T|^2\big) =  f^\tp \Sigma_T f
\]
\revblue{Because of this interpretation of the minimum m.s.e. as the conditional variance, the estimator $f^\tp m_T$ is referred to as the minimum variance estimator.  Our interest is to compute it from solving an optimal control problem.  For this purpose, one considers an estimator that is a weighted linear combination of past observations. Such an estimator is referred to as linear predictor.  In continuous-time settings, a linear predictor is of the form
\[
S_T =  \text{(constant)} - \int_0^T u_t^\tp \ud Z_t
\]
where the control input, denoted by $u:= \{u_t\in\Re^m:0\le t \le T\}$, determines the weights used to combine past observations.  The control input $u$ is an element of the function space $L^2([0,T];\Re^m)$ and parametrizes the estimator: each choice of $u$ yields a distinct $S_T$.  The problem is to choose an optimal such $u$, provided one exists, so as to obtain the minimum variance estimator.  For this purpose, it is mathematically convenient to introduce a dual process $\{y_t\in\Re^d:0\leq t\leq T\}$. As shown in the sidebar ``Duality principle calculations'', the dual process allows one to express the m.s.e. as an optimal control-type objective.  The mathematical form of the resulting optimal control problem is as follows: }

\begin{itemize}
\item Minimum variance dual optimal control problem:
\end{itemize}
\begin{subequations}\label{eq:mv-intro}
	\begin{align}
		\mathop{\text{Minimize}}_{u\in L^2([0,T];\Re^m)}\!:\quad \sJ_T(u) &= |y_0|^2_{\Sigma_0}+ \int_0^{T} |y_t|^2_Q + |u_t|_R^2 \ud t \label{eq:mv-intro-a}\\
		\text{Subject to}\;\;:\; -\frac{\ud y_t}{\ud t}
		&= A y_t + H u_t,\quad y_T = f \label{eq:mv-intro-b}
	\end{align}
\end{subequations}

\revblue{Note the terminal condition ``$y_T = f$'' in~\eqref{eq:mv-intro-b} is according to the given vector $f\in\Re^d$ for which $\E(f^\tp X_T\mid \clZ_T)$ is being sought.  Based on the calculations in the sidebar, the relationship between optimal filtering and optimal control problems is formally given by the following duality principle (see~\cite[Sec.~3.4.3]{JinPhDthesis} for a complete proof):}

\begin{theorem}[Duality principle for the linear-Gaussian model~\eqref{eq:linear-Gaussian-model}] \label{prop:duality-KalmanBucy}
	For a control input $u \in L^2([0,T];\Re^m)$, consider an estimator
	\begin{equation}\label{eq:LG-estimator}
		S_T:= y_0^\tp m_0 - \int_0^T u_t^\tp \ud Z_t
	\end{equation}
	Then
	\begin{equation*}\label{eq:LG-duality-principle}
		\bsJ_T(u) = \E\big(|f^\tp X_T - S_T|^2\big) 
	\end{equation*}
Suppose $u^\opt=\{u_t^\opt\in\Re^m:0\le t \le T\}$ is 
the optimal control input and
$\{y_t:0\leq t\leq T\}$ is the resulting optimal trajectory.  Then 
\begin{equation}
y_t^\tp m_t = y_0^\tp m_0 - \int_0^t (u_s^\opt)^\tp \ud
Z_s,\quad 0\leq t\leq T\label{eq:mean_KF_optimal_control}
\end{equation}
and therefore, at the terminal time $t=T$, $f^\tp m_T = y_0^\tp m_0 - \int_0^T (u_s^\opt)^\tp \ud
Z_s$ with the optimal value $\bsJ_T(u^\opt) =  f^\tp \Sigma_T f$.
\end{theorem}
The sidebar contains a derivation (following~\cite[Sec.~7.6]{astrom1970}) of the equation of the Kalman filter~\eqref{eq:KF_eqn} starting from the solution of the minimum variance 
optimal control problem~\eqref{eq:mv-intro}. \revblue{While minimum variance duality is classical, a modern interpretation in terms of reproducing kernel Hilbert space appears in~\cite{aubin2022reproducing}.}

\subsection{The issue with reversing the arrow of time}

We have now reviewed the two manners to express the
dual relationship between estimation and control:
\begin{romannum}
\item Elementary duality between observability of the state-output
  system~\eqref{eq:LTI-obs} and the controllability of the input-state system~\eqref{eq:LTI-ctrl}; and
\item Minimum variance duality that relates the optimal filtering
  problem for the linear Gaussian
  model~\eqref{eq:linear-Gaussian-model} to a deterministic LQ
  optimal control problem for the control system~\eqref{eq:mv-intro-b}.
\end{romannum}
We begin by commenting upon similarities and differences between the
elementary duality and the minimum variance duality: An obvious
similarity is the appearance of the same backward-in-time dual control
system,~\eqref{eq:LTI-ctrl} and~\eqref{eq:mv-intro-b}, respectively, for the two
cases.  An important difference is that while the output $z$
of~\eqref{eq:LTI-obs} is a
deterministic function of time the observation process
$Z$ is a stochastic process.  The appropriate function space for
deterministic $z$ is $L^2([0,T];\Re^m)$.  Because $Z$ is a stochastic process,
the appropriate function space is $L^2_\clZ([0,T];\Re^m)$: the
  function space of $\Re^m$-valued stochastic processes on $[0,T]$ that are
  forward-adapted to the filtration $\clZ$.
  Formally, an element of this function space is any $\Re^m$-valued stochastic
  process $U=\{U_t:0\leq t\leq T\}$ whose value $U_t$ at time $t$ is completely determined
  from past observations up to time $t$.  A rigorous definition
  appears at a later point in this paper (see also~\cite[Sec.~4.1]{JinPhDthesis}).
A practical
consequence of this difference is that one may then consider more general
estimators of the form (compare with~\eqref{eq:LG-estimator})
	\begin{equation}\label{eq:LG-estimator-stoch}
		S_T:= \text{(constant)} - \int_0^T U_t^\tp \ud Z_t
	\end{equation}
where $U\in L^2_\clZ([0,T];\Re^m)$ (the righthand-side is well defined as an It\^{o} stochastic integral).  
Because the space of deterministic control inputs is a subspace of
this larger function space, it may even be possible to obtain a
solution with smaller m.s.e. by considering such estimators (or
at least the possibility can not be ruled out a priori)!

\begin{figure}
	\centering
	\includegraphics[width=0.9\textwidth]{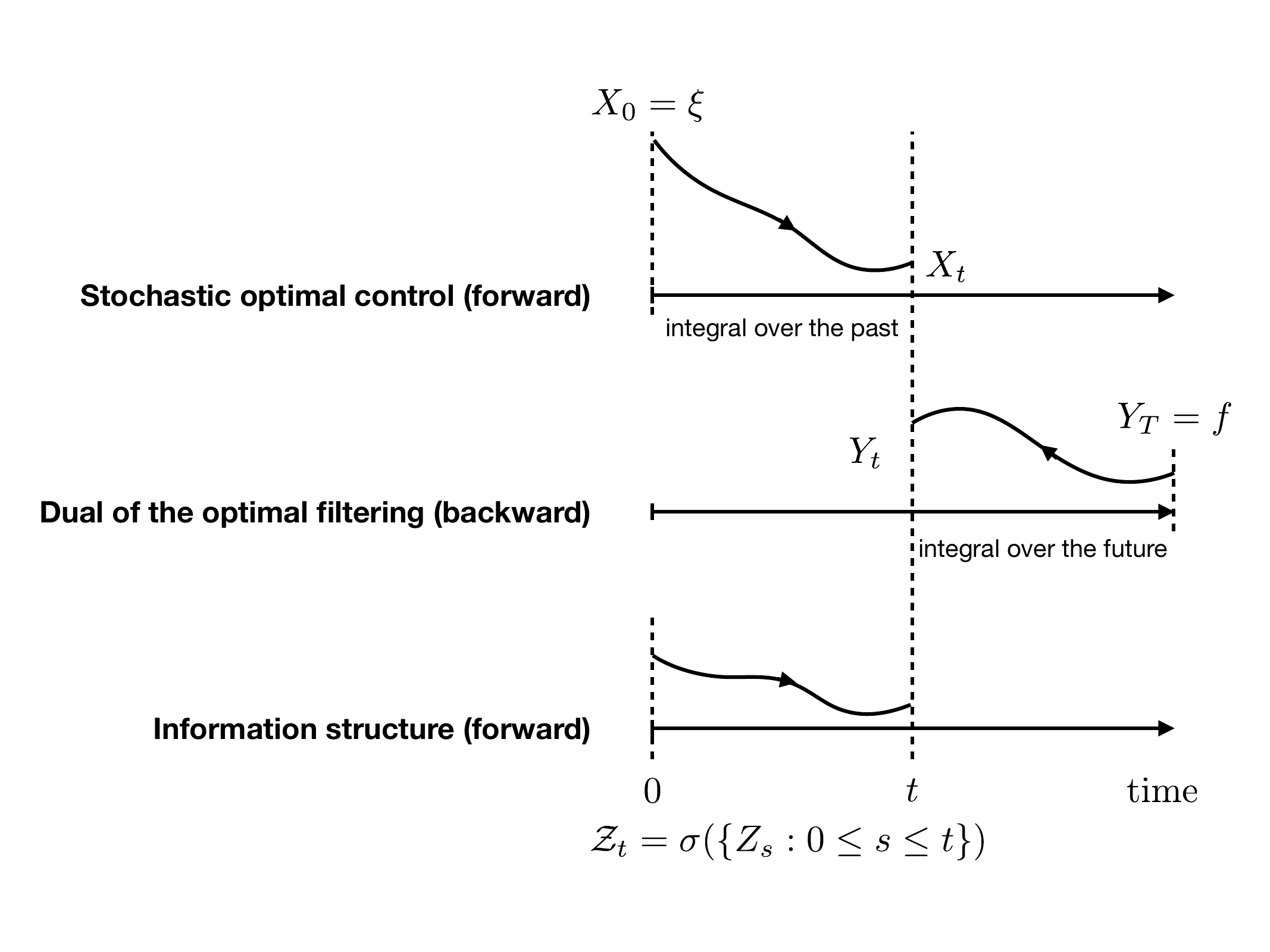}
	\caption{Arrow of time: (top) in stochastic
          optimal control, a control system is integrated
          forward-in-time from initial condition $X_0=x$; (middle)
          dual optimal control system is integrated backward-in-time
          from terminal condition $Y_T=f$.  In either case, the
          control input $U_t$ at time $t$ is allowed to depend only
          upon past observations -- captured using (bottom) the
          information structure (filtration) which evolves
          forward-in-time.  This is key difficulty in extending
          Kalman's duality to nonlinear non-Gaussian settings. 
} \label{fig:information-structure}
\end{figure}

A difficulty that now arises is that~\eqref{eq:mv-intro-b}
is a deterministic ODE.  With a stochastic $U$, its solution will
also be a stochastic process.  For this reason, let us consider a
naive extension of the dual control system as follows:
\begin{equation}\label{eq:naiveY}
-\ud Y_t = (A Y_t + H U_t)\ud t,\quad Y_T=f
\end{equation}
where the capital letter $Y=\{Y_t:0\leq
t\leq T\}$ is used to denote the fact that the dual process is now a stochastic
process. However, this poses a problem of causality, namely, the
righthand-side $A Y_t + H U_t$ is a sum of two terms: $AY_t$
which depends upon the future (because of the backward-in-time nature
of~\eqref{eq:naiveY}), and $H U_t$ which depends upon the past (because
it is an element of $L^2_\clZ([0,T];\Re^m)$). \revblue{The problem is depicted in Figure~\ref{fig:information-structure} based on drawing a contrast with stochastic control.}

An intuitive explanation for this difficulty, which also suggests a way
to its resolution, is that the apparent
similarity between the deterministic state-output
system~\eqref{eq:LTI-obs} and the stochastic linear Gaussian
model~\eqref{eq:linear-Gaussian-model} is deceptive.  Because of the
time reversal, it is not at all clear that~\eqref{eq:LTI-ctrl} is {\em also} the appropriate dual control
system for the stochastic model~\eqref{eq:linear-Gaussian-model}.
As it turns out, for the linear Gaussian model, it suffices to
consider deterministic control inputs.  Therefore, the issue of
time reversal becomes moot and the minimum
variance duality was already well understood
in 1960s. 
For the more general HMM, the issue of time reversal needs to be
carefully addressed.  The resolution appears in the PhD dissertation
of the first author~\cite{JinPhDthesis} and is described in this article.


\revblue{
Before presenting the generalization, we have one more task which is to describe the minimum energy duality.  At the onset, we noted that duality is expressed in two inter-related manners:
\begin{romannum}
\item Duality between observability and controllability; and
\item Duality between optimal filtering and optimal control.
\end{romannum}
Elementary duality between~\eqref{eq:LTI-obs} and~\eqref{eq:LTI-ctrl} is an example of the former (item (i) in the list above) while the minimum variance duality~\eqref{eq:mv-intro} is an example of latter (item (ii) in the list). As it turns out, minimum variance is not the {\em only} way to express an optimal filtering problem as an optimal control problem.  Another way to do so is the minimum energy duality which is described next.}  The sidebar ``Simple example to illustrate the two types of duality'' is included to explain and relate the two types of duality for the simplest example of estimating a Gaussian random variable.

\subsection{Minimum energy duality for the linear Gaussian model}

Like the minimum variance duality, the classical form of the minimum energy duality too is for the linear Gaussian model~\eqref{eq:linear-Gaussian-model} over a finite time-horizon $[0,T]$.  While the minimum variance duality starts with the consideration of the m.s.e., the objective of interest in the minimum energy duality is the (minus of the) {\em log-likelihood function} defined as follows:
\[
-2\log \rho_{X\mid Z}(x\mid z)
\]
where $\rho_{X\mid Z}$ is the conditional density for the linear
Gaussian model~\eqref{eq:linear-Gaussian-model} and $x = \{x_t:0\le t
\le T\}$ and $z = \{z_t:0\le t \le T\}$ are the possible state and
output trajectories.  Based on a calculation in the 1963 paper of
Bryson and Frazier and a more
widely known (in the duality literature) 1968 paper of
Mortensen,
the objective is expressed as a quadratic optimal control-type
cost~\cite{bryson1963smoothing,mortensen1968}.  The control input $v:=\{v_t\in\Re^n:0\leq t\leq T\}$
now has the dimension $n$ of the process noise \revblue{(which is different from the dimension $m$ of the control input $u$ in the minimum variance duality)}. The dual optimal control
problem is as follows:

\begin{itemize}
\item Minimum energy dual optimal control problem:
\end{itemize} 
\begin{subequations}\label{eq:mee-intro}
	\begin{align}
		\mathop{\text{Minimize}}_{\stackrel{v\in L^2([0,T];\Re^n)}{x_0\in\Re^d}}\!:\quad &
		\sJ_T(v,x_0;z)=
		|m_0-x_0|^2_{\Sigma_0^{-1}}
	+ \int_0^{T} |v_t|^2+ |\dot{z}_t - H^\tp x_t|_{R^{-1}}^2 \ud t \label{eq:mee-intro-a}\\
		\text{Subject to}\;\;:\;& \frac{\ud x_t}{\ud t}
		= A^\tp x_t + \sigma v_t \label{eq:mee-intro-b}
	\end{align}
\end{subequations}
where $z = \{z_t\in\Re^m:0\leq t\leq T\}$ is a given sample path of
observation data.   


The optimal solution is 
given in the following proposition (see~\cite[Sec.~3.4.5]{JinPhDthesis} for a proof):

\begin{proposition}\label{prop:min-energy-KalmanBucy}
	Consider the optimal control problem~\eqref{eq:mee-intro}.
	Then for any choice of $v$ and $x_0$, 
	\[
	\bsJ_T(v,x_0;\dot{z}) \ge \int_0^T |\dot{z}_t-H^\tp \hat{x}_t|_{R^{-1}}^2 \ud t
	\]
	where the process $\hat{x} = \{\hat{x}_t:0\le t\le T\}$ is
        given by
\begin{subequations}\label{eq:LG-smoother}
	\begin{equation}\label{eq:LG-smoother-forward}
		\frac{\ud \hat{x}_t}{\ud t} = A^\tp \hat{x}_t +
                \Sigma_tH R^{-1} (\dot{z}_t - H^\tp \hat{x}_t) ,\quad 0\leq t\leq T, \qquad \hat{x}_0 = m_0
	\end{equation}
(Compare with Kalman filter~\eqref{eq:KF_mean} for the conditional mean
$\{m_t:t\geq 0\}$). 
Assuming $\Sigma_t$ is invertible for all $0\leq t\leq T$, 
	the equality holds with the optimal trajectory given by the backward-in-time equation:
	\begin{equation}\label{eq:LG-smoother-backward}
		\frac{\ud x_t^\opt}{\ud t} = A^\tp x_t^\opt + \sigma \underbrace{\sigma^\tp \Sigma_t^{-1}(x_t^\opt-\hat{x}_t)}_{v_t^\opt},\quad 0\leq t\leq T, \qquad x_T^\opt = \hat{x}_T
	\end{equation}
\end{subequations}
\end{proposition}

\revblue{
  The minimum energy duality~\eqref{eq:mee-intro} is by far the more common manner to express an optimal estimation problem as an optimal control problem.  Even though the objective~\eqref{eq:mee-intro-a} is derived starting from the stochastic model~\eqref{eq:linear-Gaussian-model}, it is appropriate to think of~\eqref{eq:mee-intro} as a deterministic estimation problem.  The reasons are explained as part of the sidebar ``Minimum energy duality as a deterministic optimal control problem'', where some historical context is also provided.
  }

\subsection{Comparison of the minimum variance and minimum energy duality}  

As evident from the construction, the two types of duality
correspond to distinct estimation objectives: 
\begin{itemize}
\item Minimum variance~\eqref{eq:mv-intro} is the appropriate duality for the filtering
problem: Formula~\eqref{eq:mean_KF_optimal_control} shows that a
portion of the 
optimal control $\{u_s^{\text{OPT}}:0\leq s\leq t\}$ yields the filter
at time $t$.  Although not explicitly given as part of \Prop{prop:duality-KalmanBucy}, the optimal
cost-to-go at time $t$ has the interpretation of the optimal
m.s.e..
\item Minimum energy~\eqref{eq:mee-intro} is the appropriate duality for the
  smoothing problem: This is because the optimal trajectory $x_t =
  \E(X_t|\clZ_T)$ for $0\leq t\leq T$ (equality is meant 
  in a pathwise sense for each fixed sample-path $Z_t(\omega)=z_t$).   
The forward-backward equation~\eqref{eq:LG-smoother} is well known
(see~\cite[Eq.~(16.5.11)]{kailath2000linear}) as
the Rauch-Tung-Streibel (RTS) smoother~\cite{rauch1965maximum}.  Notably, the optimal solution at time $t$ depends
upon both the past and future values of the output.  
\end{itemize}

Let us next compare these 
vis-\`a-vis
Kalman's PoD:
\begin{romannum}
\item In the minimum variance duality, the
  constraint~\eqref{eq:mv-intro-b} is the dual control system (with
  matrix transposes).  In
  contrast, the constraint~\eqref{eq:mee-intro-b} in the minimum
  energy duality is a modified copy of the 
  model~\eqref{eq:linear-Gaussian-model-a}.  
\item In~\eqref{eq:mv-intro}, the control input $u$ has the same
  dimension as the output process while in~\eqref{eq:mee-intro},
  the control input \revblue{$v$} has the dimension of the process noise.
  Note that in the elementary duality between controllability and
  observability, the inputs and outputs are the elements of the same
  function space $L^2([0,T];\Re^m)$ and as such have the same dimension.  
\item The constraint~\eqref{eq:mv-intro-b} is
  a backward-in-time ODE while the constraint~\eqref{eq:mee-intro-b}
  is a modified copy of the signal model which proceeds forward-in-time.  
\end{romannum}
All of this shows that the minimum variance duality is 
consistent with each of the three steps in the Kalman's PoD.  In contrast, the
minimum energy duality is not consistent with any of the three steps.  The practical import of this is as
follows: The condition for
  asymptotic analysis of the minimum variance dual optimal control problem~\eqref{eq:mv-intro} is stabilizability
  of~\eqref{eq:mv-intro-b}, and by duality detectability of
  $(A^\tp,H^\tp)$.  
The latter is well known to be the appropriate condition 
  for stability of the Kalman filter~\cite{kwakernaak1969linear,ocone1996asymptotic}. In contrast,
  for minimum energy dual optimal control problem~\eqref{eq:mee-intro}, the
  important 
  condition again is detectability of $(A^\tp,H^\tp)$ but that is not at
  all 
  transparent from the constraint~\eqref{eq:mee-intro-b}. With $\sigma=0$, which is a perfectly reasonable choice for the linear
  Gaussian model~\eqref{eq:linear-Gaussian-model}, the constraint~\eqref{eq:mee-intro-b} is not even a control system.


Because the constraint~\eqref{eq:mee-intro-b} is a modified copy of the model~\eqref{eq:linear-Gaussian-model-a}, there have been several extensions of the minimum energy duality to both deterministic and stochastic nonlinear systems.  These are reviewed and compared with our original work on the generalization of minimum variance duality described next.

\section{Duality for HMMs with white noise observation model}

\revblue{
In this section, our goal is to introduce a new generalization of the linear Gaussian minimum variance duality, specifically extending the dual optimal control problem~\eqref{eq:mv-intro} and its solution outlined in~\Thm{prop:duality-KalmanBucy} to a class of nonlinear, non-Gaussian models (HMMs). The presentation closely mimics the presentation in the prior section: the generalization of~\eqref{eq:mv-intro} and~\Thm{prop:duality-KalmanBucy} are presented in~\eqref{eq:dual-optimal-control} and~\Thm{thm:duality-principle-HMM}, respectively. We begin with a description of the model.} 

\subsection{Hidden Markov model (HMM)}

Consider a pair of continuous-time stochastic processes $(X,Z)$
as follows:
\begin{itemize}
	\item The \emph{state process} $X = \{X_t:0\le t \le
          T\}$ is a Feller-Markov process taking values in the state-space
          $\bS$.  The time-evolution of a Markov process is modeled
          through its infinitesimal generator, denoted by $\clA$, and
          a prior, denoted by $\mu$.  
          The prior $\mu$ is an element of the space of
          probability measures on $\bS$.  This space is denoted by
          $\clP(\bS)$. 
\item  The \emph{observation process} $Z = \{Z_t:0\le t \le T\}$ is
  $\Re^m$-valued and satisfies the SDE (compare with~\eqref{eq:linear-Gaussian-model-b})
	\begin{equation}\label{eq:obs-model}
		\ud Z_t = h(X_t) \ud t + \ud W_t,\quad Z_0=0
	\end{equation}
	where $h:\bS\to \Re^m$ is referred to as the observation function and $W =
        \{W_t:0\le t \le T\}$ is an $m$-dimensional B.M. with covariance matrix $R
                                 \succ 0$.
	It is assumed that $W$ is independent of $X$.
\end{itemize}
The above is referred to as an HMM with \emph{white noise observations}. The model is denoted by 
$(\clA,h)$.  In applications, the important examples of HMM are as
follows:
\begin{romannum}
\item HMM on a finite state-space $\bS = \{1,2,\hdots,d\}$ of
  cardinality $d$. \revblue{For the finite state-space model, the generator $\clA$ is identified with a $d\times d$ transition rate matrix denoted as $A$.}
\item Nonlinear SDE or an ODE on the state-space
  $\bS\subseteq\Re^d$. A SDE model is referred to as an {\em It\^o
    diffusion}. The linear Gaussian
  model~\eqref{eq:linear-Gaussian-model} is an example of an It\^o
  diffusion on $\bS=\Re^d$.
\end{romannum}
The notation for these models is described in a self-contained manner in the accompanying sidebar ``Examples of HMM''.
  When $\bS$ is not finite,
  additional assumptions on the state-space $\bS$ and the model $(\clA,h)$ are necessary for well-posedness reasons. 

In the following two subsections, our goal is as follows:
\begin{romannum}
\item Define a dual control system for the HMM  $(\clA,h)$.
\item Describe the minimum variance duality and the duality principle
  for the same.
\end{romannum}
After these have been described, the classical linear Gaussian duality
is shown to be arise as the special case in the final subsection.

\subsection{Dual control system}

We begin by revisiting the sidebar ``Mathematical foundations of
duality''.
The first step is to define the function spaces.  These
should already be familiar to the reader for the state-output
linear system~\eqref{eq:LTI-obs}.  The counterparts of these for the
HMM are as follows:
\begin{romannum}
\item \textbf{Function space for the state (initial condition)}.  The initial
  condition for an HMM is the prior $\mu$.  The prior is an element of
  $\clP(\bS)$ which is not a vector space.  For functional
  analytic reasons, one regards $\clP(\bS)$ as a subset of a measure
  space denoted as $\clY^\dagger$. If $\bS$ is finite of cardinality
  $d$ then $\clY^\dagger=\Re^d$ and $\clP(\bS)$ is the probability simplex in $\Re^d$. 
\item \textbf{Function space for the signal (output)}.  For an HMM, the observation process
  $Z$ is a stochastic process.  It is an element of the Hilbert space
  $L^2_\clZ([0,T];\Re^m)$ whose mathematical definition appears in the sidebar ``Duality between BSDE control
  system and HMM''. 
\end{romannum}

For the dual control system, the function spaces are
obtained through duality pairing. 
\begin{romannum}
\item Because the state for an HMM is a measure on $\bS$,
  the state for the dual control system is a function on $\bS$.
  The space of such functions is denoted as $\clY$ with $\clY^\dagger$
  as the (dual) space of linear functionals on $\clY$. If $\bS$ is
  finite of cardinality $d$ then $\clY=\Re^d$.
\item Because $L^2_\clZ([0,T];\Re^m)$ is a Hilbert space, the control
  inputs for the dual control system are also elements of
  $L^2_\clZ([0,T];\Re^m)=:\clU$. 
\end{romannum}

The main result
of~\cite{duality_jrnl_paper_I} is to introduce the following linear backward
stochastic differential equation (BSDE) as the dual control system for
the HMM $(\clA,h)$:
\begin{subequations}\label{eq:dual-bsde}
	\begin{align}
		-\ud Y_t(x) &= \big((\clA Y_t)(x) +
                              h^\tp(x)(U_t+V_t(x))\big)\ud t - V_t^\tp
                              (x) \ud Z_t,\quad x\in \bS, \quad 0\leq t\leq T \label{eq:dual-bsde-a}\\
		Y_T(x) &= c,\quad\forall\,x\in\bS\label{eq:dual-bsde-b}
	\end{align}
\end{subequations}
where the control input $U:=\{U_t: 0\leq t\leq T\}\in \clU$ and $c\in
\Re$ is a deterministic constant.  The solution of the BSDE $(Y,V)
:=\{(Y_t(x),V_t(x)) : 0\leq t\leq T,\;x\in\bS\}$ is a $\clY\times
\clY^m$-valued stochastic process. The process is 
forward-adapted to the filtration $\clZ$.  In the settings where
the state-space $\bS=\Re^d$,~\eqref{eq:dual-bsde} is a backward stochastic
partial differential equation (BSPDE).  The well-posedness results for finite
state-space can be found in~\cite[Ch.~7]{yong1999stochastic} and for
the BSPDE in~\cite{ma1997adapted} (see also~\cite[Rem.~7-9]{duality_jrnl_paper_I}). \revblue{While the well-posedness theory is well developed, it is generally much more difficult to numerically solve a BSDE compared to a (forward) SDE. Such is the case even for the finite state-space where $\clY=\Re^d$.  The problem is that, outside some special cases, there is no explicit formula for the process $V$~\cite{pardoux2014stochastic}.}

In the  sidebar ``Duality between BSDE control system and HMM'', the mathematical characterization of duality is presented. Based on this characterization, Table~\ref{tb:comparison0} provides a side-by-side comparison showing the parallels between the classical duality for deterministic linear systems (\eqref{eq:LTI-obs} and~\eqref{eq:LTI-ctrl}) and the duality for the HMM $(\clA,h)$.  Additional
details can be found in the journal paper~\cite{duality_jrnl_paper_I}.


Our goal now is
to use the dual control system to construct the minimum variance
optimal control problem.

\begin{table*}[t]
	\centering
	\renewcommand{\arraystretch}{1.7}
	\caption{Comparison of the controllability--observability duality for linear and nonlinear systems} \label{tb:comparison}
	\small
	\begin{tabular}{p{0.14\textwidth}p{0.37\textwidth}p{0.4\textwidth}} \label{tb:comparison0}
		& {\bf Linear deterministic systems} & {\bf Nonlinear stochastic systems} \\ \hline \hline
		Function space for inputs and outputs & $\clU =  L^2([0,T];\Re^m)$\vspace{4pt} \newline  $\langle u,v\rangle = \displaystyle\int_0^T u_t^\tp v_t\ud t$ \vspace{3pt} & $\clU =  L^2_\clZ(\Omega\times [0,T];\Re^m)
		$\vspace{4pt} \newline $\langle U,V\rangle = \displaystyle\hE\Big(\int_0^T U_t^\tp V_t\ud t\Big)$\vspace{3pt} \\ \hline
		Function space for the dual state & $\clY =\Re^d$ \vspace{3pt} \newline $\langle x,y \rangle = x^\tp y$ & $\clY = C_b(\bS)$, $\clY^\dagger = {\cal M}(\bS)$ \vspace{3pt} \newline $\langle \mu,y\rangle = \mu(y)$ \\ \hline
		Linear operators & $\clL : \clU \to \clY$, $u\mapsto y_0$ by ODE~\eqref{eq:LTI-ctrl} \vspace{2pt} \newline 
		$\clL^\dagger : \clY \to \clU$, 
                  \vspace{2pt} \newline \phantom{$\clL^\dagger$}
                  \hspace{0.15em} $x_0\mapsto
                                  \{z_t:0\leq t\leq T\}$ by ODE~\eqref{eq:LTI-obs}
		& $\clL:\clU\times\Re\to\clY$, $(U,c)\mapsto Y_0$ by BSDE~\eqref{eq:dual-bsde} \vspace{2pt}\newline
		$\clL^\dagger:\clY^\dagger\to \clU\times\Re$,
                  \vspace{2pt} \newline \phantom{$\clL^\dagger$}
                  \hspace{0.15em} $\mu\mapsto
                  (\{\sigma_t(h):0\leq t\leq T\}, \mu(\ones))$
                  \vspace{2pt} \newline \phantom{$\clL^\dagger$}
                  \hspace{0.15em} by the Zakai equation~\eqref{eq:Zakai}
		\\ \hline
		Controllability & $\Rsp(\clL) = \Re^d$ & $\overline{\Rsp(\clL)} = \clY$ \\ \hline
		Observability & $\Nsp(\clL^\dagger) = \{0\}$ & $\Nsp(\clL^\dagger) = \{0\}$ \\ \hline
		Duality & \multicolumn{2}{l}{$\Rsp(\clL)^\bot =
                          \Nsp(\clL^\dagger) \quad \Longrightarrow
                          \quad$ A state-output system is observable iff
                          the dual control system is controllable.
	}
		\\ \hline \hline
	\end{tabular}
\end{table*}




\subsection{Minimum variance duality for HMM}

Recall that for the linear Gaussian
model~\eqref{eq:linear-Gaussian-model}, the filtering objective is to
compute the conditional expectation $\E(f^\tp X_T|\clZ_T)$.
The function $f^\tp x$ is an example of a linear function of the
state $x\in\Re^d$.  In nonlinear (or stochastic) filtering, the
objective is to compute an estimate of $f(X_T)$ where the function $f\in
\clY$.  The theory allows for consideration of also random functions
that are $\clZ_T$ measurable.  For any such random
function $F$, the nonlinear filter is defined as a conditional
expectation as follows:
\[
\text{(nonlinear filter)} \qquad \pi_T(F) := \E (F(X_T)|\clZ_T),\quad F\in\clZ_T
\]
where ``$F\in\clZ_T$'' means that the function $F$ is allowed to be 
$\clZ_T$-measurable (i.e., for each $x\in\bS$, $F_T(x)$ is allowed to depend upon
observations up to time $T$). As before, the conditional expectation has the
following interpretation as the solution of an optimization problem:
\begin{equation*}\label{eq:minimum-variance}
\pi_T(F) = \mathop{\operatorname{argmin}}_{S_T\in \clZ_T} \E\big(|F(X_T)-S_T|^2\big)
\end{equation*}
The dual control system is used to express the above as a dual optimal control problem.

\begin{itemize}
\item Minimum variance optimal control problem:
\end{itemize}
\begin{subequations}\label{eq:dual-optimal-control}
	\begin{align}
&           \mathop{\text{Minimize:}}_{U\in\;\clU}\; \sJ_T(U)  =    \E (|Y_0(X_0)-\mu(Y_0)|^2) + 
\E \Big(\int_0^T l (Y_t,V_t,U_t\,;X_t) \ud t \Big)
\label{eq:dual-optimal-control-a}
\\
&          \text{Subject to (BSDE constraint):} \nonumber \\ 
&  -\!\ud Y_t(x) = \big((\clA Y_t)(x) + h^\tp (x) (U_t +
            V_t(x))\big)\ud t - V_t^\tp(x)\ud Z_t,\;\; x\in\bS,\;\; 0\leq t\leq T \quad \\ & \quad\;\; Y_T (x)  = F(x), \;\; x \in \bS 
\label{eq:dual-optimal-control-b}
	\end{align}
\end{subequations}
where the running cost 
\begin{equation*}\label{eq:running_cost_formula}
l(y,v,u;x):= (\Gamma y)(x) + |u+v(x)|_R^2,\quad y\in\clY,\;v\in\clY^m,\;u\in\Re^m,\;x\in\bS
\end{equation*}
Here, $\Gamma$ is the \emph{carr\'e du
  champ} operator whose definition and meaning is described in the
sidebar ``Carr\'e du champ operator''. \eqref{eq:dual-optimal-control}
is an example of a linear quadratic stochastic optimal control problem
subject to a BSDE constraint.  If the state-space is finite, the
problem is finite-dimensional while for $\bS=\Re^d$, the problem
is infinite-dimensional.  In either case, its solution based on the maximum principle appears
in~\cite[Thm.~2]{duality_jrnl_paper_II}.

The relationship between optimal filtering and optimal control
problems is given by the following duality principle (compare
with~\Thm{prop:duality-KalmanBucy}) which is an original contribution
of our work: 

\medskip

\begin{theorem}[Duality principle for the HMM $(\clA,h)$~\cite{duality_jrnl_paper_II}]\label{thm:duality-principle-HMM}
	For a control input $U\in \clU$, consider an estimator
	\begin{equation}\label{eq:estimator}
		S_T := \mu(Y_0) - \int_0^T U_t^\tp \ud Z_t
	\end{equation}
	Then 
	\begin{equation}\label{eq:duality-principle}
		\sJ_T(U) = \E\big(|F(X_T)-S_T|^2\big)
	\end{equation}
Suppose $U^\opt=\{U_t^\opt:0\le t \le T\}$ is
the optimal control input and
$\{Y_t:0\leq t\leq T\}$ is the resulting optimal trajectory.  Then 
\[
\pi_t(Y_t) = \mu(Y_0)  - \int_0^t (U_s^\opt)^\tp \ud
Z_s,\quad 0\leq t\leq T
\]
and therefore, at the terminal time $t=T$, 
\[
\pi_T(F) = \mu(Y_0)  - \int_0^T (U_s^\opt)^\tp \ud
Z_s
\]
with the optimal value $\bsJ_T(U^\opt) =  \E\big(|F(X_T)-\pi_T(F)|^2\big)$.
\end{theorem}

The duality principle first appeared in our paper~\cite[Thm.~1 and Prop.~1]{kim2019duality}.  \revblue{A proof sketch for~\eqref{eq:duality-principle} appears as part of the sidebar.} 
As in the linear Gaussian case, minimum
variance optimal control problem~\eqref{eq:dual-optimal-control} is the appropriate duality for the 
filtering problem.
The optimal control input not only yields the filter at the terminal
time $T$ but also a portion of the optimal control $\{U_s^\opt:0\leq s\leq
t\}$ yields the filter at time $t$.

A special case
of~\eqref{eq:dual-optimal-control} based on deterministic control
input $U=u$ and terminal condition $F=f$ is considered in the
M.S. thesis~\cite{jan2021master}.  With such a choice, $V=0$ and the
analysis is simpler but the estimator is
sub-optimal~\cite[Rem.~1]{kim2019duality}. An equivalent formulation of~\eqref{eq:dual-optimal-control} expressed as a forward-backward SDE appears in~\cite{kim2024FBSDE}.   
	
In the following, the classical minimum variance duality~\eqref{eq:mv-intro} for the linear Gaussian model is shown to be a special case of this more general result.

\subsection{Linear Gaussian special case}

Consider the linear Gaussian model~\eqref{eq:linear-Gaussian-model}.
We impose the following restrictions:
	\begin{itemize}
		\item The control input $U=u$ is restricted to be a deterministic
		function of time.  In particular, it does not depend upon the
		observations.  
		For such a control input, the solution $Y=y$ of the
                BSPDE is a deterministic function of time, and $V=0$.
                The BSPDE becomes a PDE:
		\begin{equation}\label{eq:det_pde}
			-\frac{\partial y_t}{\partial t}(x) = (\clA
                        y_t)(x) + x^\tp H 
			u_t,\quad y_T(x) = c,\quad x\in \Re^d
		\end{equation}
		where the lower-case notation is used to stress the fact that $u$
		and $y$ are now deterministic functions of time.
		
		\item Consider a
                  finite ($d$-) dimensional space of linear functions:
		\[
		{\cal S}:=\{f:\Re^d\to \R \;:\;  f(x)= \tilde{f}^\tp x\;\text{where}\; \tilde{f}\in\Re^d\}
		\]  
		Then ${\cal S}$ is an invariant subspace for the
		dynamics~\eqref{eq:det_pde}.  On ${\cal S}$, expressing
                $y_t(x) = \tilde{y}_t^\tp x$, using the formula~\eqref{eq:gen-LG}, the PDE reduces to an ODE:  
		\begin{equation}\label{eq:LTI-ctrl_1}
		-\frac{\ud \tilde{y}_t}{\ud t}
		= A \tilde{y}_t + H u_t,\quad \tilde{y}_T = 0
		\end{equation}
                where the terminal condition is set to 0 because it is the
                only constant function which is also linear.
	\end{itemize}
In this manner, we have recovered the dual control
system first introduced as~\eqref{eq:LTI-ctrl} in the discussion of
elementary duality between controllability and observability. 

The {\em only} assumption in going from
BSPDE~\eqref{eq:dual-bsde} to the ODE~\eqref{eq:LTI-ctrl_1} is that the
control $U$ is deterministic.  That a deterministic control suffices for minimum variance
duality is because of a standard result in the theory of Gaussian
processes that conditional expectation $\E(f^\tp  X_T|\clZ_T)$ is evaluated in the form of
a linear predictor~\cite[Cor.~1.10]{le2016brownian}.  For this reason,
it suffices to consider an estimator of the form 
\[
	S_T := \text{(constant)} - \int_0^T u_t^\tp \ud Z_t
\]
where $u = \{u_t \in \Re^m: 0\le t\le T\}$ is 
deterministic.
Consequently, for linear Gaussian estimation, it suffices to restrict the
admissible space of control inputs to $L^2\big([0,T];\Re^m\big)$ which
is a subspace of $L_{\clZ}^2\big([0,T];\Re^m\big)$.  
Using a deterministic control $u$, and the terminal condition $F(x) = f^\tp x$, the solution of the BSPDE~\eqref{eq:dual-optimal-control-b} is given by
\[
Y_t(x) = y_t^\tp x, \quad V_t(x)=0, \quad x\in \Re^d, \;\;0\leq t\leq T
\] 
where $y=\{y_t\in\Re^d:0\leq t\leq T\}$ is a solution of the
backward-in-time 
ODE:
		\begin{equation*}
		-\frac{\ud y_t}{\ud t}
		= A y_t + H u_t,\quad y_T = f
		\end{equation*}
Using the formula~\eqref{eq:Gamma-LG} for the carr\'e du champ, 
the running cost
\begin{align*}
l (Y_t, V_t, U_t; X_t) &= (\Gamma Y_t)(X_t) + |U_t+V_t(X_t)|_R^2 =
                         y_t^\tp (\sigma\sigma^\tp) y_t + |u_t|_R^2 =
                         |y_t|_Q^2 + |u_t|_R^2
\end{align*} 
With the
Gaussian prior, the initial 
cost $\E (|y_0^\tp (X_0 - m_0)|^2)= y_0^\tp \Sigma_0 y_0$.

Combining all of the above, the minimum variance optimal control problem~\eqref{eq:dual-optimal-control}
reduces to its classical counterpart~\eqref{eq:mv-intro}.

\section{Discussion}

\subsection{Kalman's PoD and time reversal}

Because the proposed duality reduces to the classical minimum variance
duality, the relationship between the HMM $(\clA,h)$ and the BSDE is
consistent with the three steps of the Kalman's PoD.  We revisit these
as follows:
\begin{romannum}
\item The first step about taking transpose of matrices is not 
  applicable to the nonlinear setting.  Having said that, both the HMM and the
  BSDE are parametrized by the model $(\clA,h)$.  For the linear
  system, the model is defined by the matrices $(A,H)$.
\item The inputs and outputs constraints are interchanged.  In the
  nonlinear case, the
  constraint is that these are both elements of the common
  function space $\clU$.
\item The arrow of time is reversed. While the HMM is integrated
  forward-in-time, the BSDE is integrated backward-in-time.  
\end{romannum}

\subsection{Comparison with minimum energy duality}

\begin{table*}
	\centering
	\renewcommand{\arraystretch}{2}
	\caption{Comparison of the Mitter-Newton duality and the
          duality proposed in our work} \label{tb:comparison}
	\small
	\begin{tabular}{p{0.25\textwidth}p{0.35\textwidth}p{0.35\textwidth}}
		& {\bf Mitter-Newton duality} & {\bf Duality proposed
                                                in our work} \\
          \hline \hline
		Filtering/smoothing obj. & Minimize relative
                                                entropy (K-L)
                                              & Minimize
                                                variance (m.s.e.)
          \\ \hline
		Observation (output) process & Pathwise ($z$ is a sample path) & $Z$ is a
          stochastic process\\ \hline
		Control (input) process & $U_t$ has dimension of
                                        the process noise & $U$, 
                                                            $Z$ are
                                                            both
                                                            elem.
                                                            of $L^2_{\clZ}([0,T];\Re^m)$
                                                        \\ \hline
		Dual optimal control prob. &
                                               Eq.~\eqref{eq:opt-cont-sde-hjb-intro} 
                                              &  Eq.~\eqref{eq:dual-optimal-control}
          \\ \hline
		Arrow of time & Forward-in-time & Backward-in-time \\
          \hline  
		State-space & $\bS$: same as the state-space for $X_t$ &
                                                                      $\clY$:
                                                                             the
                                                                             space
                                                                             of functions on $\bS$ \\ \hline 
		Constraint & Controlled copy of the state
                             process
                                              & Dual control system
                                                BSDE~\eqref{eq:dual-optimal-control-b}
          \\ \hline
Running cost (Lagrangian) & $l(x,u\,;z_t)= \half |u|^2 + \half |\dot{z}_t-h(x)|_{R^{-1}}^2$ &
$l(y,v,u;x) = (\Gamma y)(x) + |u+v(x)|_R^2$ \\ \hline 
Value function (meaning) & Minus log of posterior density &
                                                                Expected value of conditional variance  \\
          \hline  
Model property  (for analysis of asymptotic stability) & Unclear & Stabilizability of BSDE \\
          \hline  
Optimal solution gives & Equation of smoothing  & 
 Equation of nonlinear filtering \\
          \hline 
		Linear-Gaussian special case & Minimum energy duality~\eqref{eq:mee-intro} & Minimum variance duality~\eqref{eq:mv-intro}  \\ \hline 
\hline
	\end{tabular}
\end{table*}

Sidebar ``Historical survey of minimum energy duality for nonlinear estimation'' contains a historical review of minimum energy duality.  The review is written to relate the early work of Mortensen and Hijab, the duality implicit in construction of the MHE algorithms, and the duality described by Mitter and Fleming and by Mitter and Newton. Table~\ref{tb:comparison} provides a side-by-side comparison of the two types of duality for nonlinear stochastic systems:
\begin{itemize}
\item Minimum energy duality of
  Mitter and Newton~\eqref{eq:opt-cont-sde-hjb-intro} on the left-hand side; and
\item Minimum variance duality~\eqref{eq:dual-optimal-control} proposed by us on the right-hand side.
\end{itemize}
The two represent generalization of the classical minimum energy duality~\eqref{eq:mee-intro} and the minimum variance
duality~\eqref{eq:mv-intro}, respectively, for the linear Gaussian
model.

The important distinctions between the two types of duality are as follows:  
\begin{romannum}
\item{\em Inputs and outputs.} The control input for the BSDE and the
  output of the HMM are dual
  processes that have the same dimension.  These are both element of
  the same Hilbert space $\clU$.  This is not true for the minimum
  energy duality.
\item{\em Stability condition.} On account of the linear quadratic nature
  of~\eqref{eq:dual-optimal-control}, stabilizability of the
  BSDE is a natural condition for asymptotic analysis.  Because 
  BSDE is the dual control system, stabilizability of the BSDE is the dual property to the
  detectability of the HMM~\cite[Cor.~2]{duality_jrnl_paper_I}.  The latter is known
  to be the appropriate condition 
  for stability of the nonlinear filter~\cite{van2010nonlinear}.
  There are no analogous conditions known for the minimum
  energy duality.
\item{\em Arrow of time.} The dual control system is backward-in-time
  consistent with Kalman's PoD.  However, note that the solution of
  the BSDE is forward-adapted to the filtration generated by the
  observations.  This is a key point concerning time reversal for the
  filtering problem.  
Adaptedness is not an issue in the minimum energy duality. 
 \end{romannum} 

\subsection{Future directions}

In this paper, we have made the case that Kalman's PoD underpins both
the elementary duality between controllability and observability of
linear systems as well as the minimum variance duality.  Both of these 
are fundamental to Control Theory and covered in the first graduate
course on the subject.  We described an original extension of these to
HMMs with white noise observations as well as the relationship of the
same with the minimum energy duality.  
A natural question that arises now is why do this?  In particular, 
what are the benefits of reformulating the optimal filtering problem
as an optimal control problem? We see two avenues for
future work that also provide an answer to this important question.

\revblue{
Expressing a dynamical system as an optimal control system is useful for the purposes of asymptotic stability analysis, as the time-horizon $T\to\infty$. 
Minimum variance duality for HMM offers an intriguing possibility to extend (Kalman filter-type) stability theory to nonlinear filters. This is a subject of continuing work from our research group at Illinois.  Results for HMM with white noise observations appear in~\cite{kim2024variance,kim2021detectable,kim2021ergodic} and the more general case is discussed in~\cite{kim2024backward}.  

The second direction is to use the minimum variance duality for the
purposes of numerical approximation of the nonlinear filter.
The nonlinear predictor form~\eqref{eq:estimator} of the estimator 
in \Thm{thm:duality-principle-HMM} is potentially applicable to a
number of application areas involving time-series prediction.  Development of algorithms in data-driven or learning type scenarios is of topical interest.
}


\newpage

\bibliographystyle{IEEEtran}
\bibliography{../../bibfiles/_master_bib_jin.bib,../../bibfiles/jin_papers.bib,../../bibfiles/extrabib.bib,../../bibfiles/estimator_controller.bib}

\sidebars 

\clearpage
\newpage
\section{Summary}
\label{sidebar:summary}

Duality between estimation and control is a foundational concept in
Control Theory.  Most students learn about the elementary duality -- 
between observability and controllability -- in their first graduate
course in linear systems theory.  Therefore, it comes as a surprise
that for a more general class of nonlinear stochastic systems (hidden
Markov models or HMMs), duality is incomplete. 

Our objective in writing this article is two-fold: (i) To describe the
difficulty in extending duality to HMMs; and (ii) To discuss its 
recent resolution by the authors.  A key message is that
the main difficulty in extending duality comes from time reversal
in going from estimation to control.  The reason for time reversal
is explained with the aid of the familiar linear deterministic and linear
Gaussian models.  The explanation is used to motivate the difference
between the linear and the nonlinear models.  Once the difference is
understood, duality for HMMs is described based on our
recent work.  The article also includes a comparison and discussion of
the different types of duality considered in literature.

\clearpage
\newpage

\section{Arrow of time in ordinary and stochastic differential equations}
\label{sidebar:arrow_of_time}

\subsection{Ordinary differential equation}

We illustrate the concepts with the aid of the simple ordinary
differential equation (ODE) over time
horizon $[0,T]$:
\begin{equation}\label{eq:ODE_arrow}
\frac{\ud y_t}{\ud t} = 0,\qquad 0 \leq t \leq T
\end{equation}
Arrow of time refers to the nature of integration.  The two types of
solutions are as follows:
\begin{subequations}\label{eq:forw_and_back_ODE}
\begin{align}
\text{(forward-in-time)} \qquad y_t &= y_0 + \int_0^t 0 \ud s,\qquad 0 \leq
                                      t \leq T  \label{eq:forw_and_back_ODE_f}\\
\text{(backward-in-time)} \qquad y_t &= y_T - \int_t^T 0 \ud s,\qquad 0
                                       \leq t \leq T 
\label{eq:forw_and_back_ODE_b}
\end{align}
\end{subequations}
It is easily verified that either of these solutions satisfy the ODE.  
The forward-in-time solution is obtained upon specifying the initial
condition $y_0=f$ and integrating forward-in-time. The
backward-in-time solution is obtained upon specifying the terminal
condition $y_T=f$ and integrating backward-in-time. Upon replacing the
right-hand side of~\eqref{eq:ODE_arrow} (which is $0$) by a suitable function of $(t,y_t)$, the
above is easily extended to the reader's favorite ODE.  

\subsection{Stochastic differential equation}

A stochastic differential equation (SDE) is different from an ODE
because its solution also depends upon a
Brownian motion (B.M.) which is denoted in this sidebar by $Z:=\{Z_t:0\leq t\leq T\}$.  In order to make the nature
of this dependence precise, one needs to introduce
$\clZ_t:=\sigma(\{Z_s:0\leq s\leq t\})$ as the sigma-algebra of
observations of $Z$ up to time $t$.  Formally, $\clZ_t$ is the set of all
events (subsets of the sample space) that become known based on 
the knowledge of B.M. up to time $t$.  The collection $\clZ:=\{\clZ_t:0\leq t\leq T\}$ is
referred to as the filtration.  A random variable $S_T$ is said to be
$\clZ_T$-measurable (written as $S_T\in \clZ_T$) if its value is known
based on the knowledge of B.M. up to time $T$. A stochastic process $Y=\{Y_t:0\leq t\leq T\}$ is
said to be forward-adapted to the filtration $\clZ$ if $Y_t\in
\clZ_t$ for $0\leq t\leq T$.  

The simplest SDE, counterpart of the ODE~\eqref{eq:ODE_arrow}, is
given by
\begin{equation}\label{eq:SDE_arrow}
\ud Y_t = 0,\qquad 0 \leq t \leq T
\end{equation}
Upon specifying the initial condition $Y_0=f$, at time $t=0$, the
forward-in-time solution is obtained as follows:
\[
\text{(forward-in-time)} \qquad Y_t = Y_0 + \int_0^t 0 \ud s + \int_0^t
0 \ud Z_s, \qquad 0 \leq t \leq  T 
\]
where ``$\int_0^t 0 \ud Z_s$'' is an example of an It\^o-integral.
The construction of the It\^o-integral ensures that the solution is
forward-adapted to $\clZ$. 

While the forward-in-time solution of the SDE is a natural generalization of
the respective formula~\eqref{eq:forw_and_back_ODE_f} for the ODE, the
backward-in-time solution is not as straightforward. The issue is that
at time $t=T$, the terminal condition $Y_T=F$ may in general depend
upon B.M. up to time $T$ (i.e., $F\in\clZ_T$).  With such a
terminal condition, a naive extension of the backward-in-time 
formula~\eqref{eq:forw_and_back_ODE_b} for the ODE gives
\begin{equation}\label{eq:wrong_backward}
Y_t = F - \int_t^T 0 \ud s - \int_t^T 0 \ud Z_s, \qquad 0 \leq t \leq T 
\end{equation}
The problem is that such a solution is not adapted to the filtration.
In particular, $Y_0=F$ depends upon future values of B.M. which are not
known at time $t=0$.

\subsection{Backward stochastic differential equation}

The theory of backward stochastic differential equation (BSDE) is
designed to define a backward-in-time solution that is
forward-adapted. The trick is to define the solution as a pair
$(Y,V)=\{(Y_t,V_t):0\leq t\leq T\}$ as follows:
\begin{equation}\label{eq:right_backward}
\text{(backward-in-time)} \qquad Y_t = F - \int_t^T 0 \ud s - \int_t^T V_s \ud Z_s, \qquad 0 \leq t \leq T 
\end{equation}
A succinct notation for the same is the BSDE
\[
-\ud Y_t = - V_t \ud Z_t, \qquad Y_T = F
\]
where the minus sign in from of $\ud Y_t$, by convention, indicates that the equation
is being integrated backward-in-time.  An important point to note is
that one solves for the pair $(Y,V)$ given the terminal condition
$Y_T=F$ at
time $T$.  The well-posedness theory of BSDE shows that a unique forward-adapted solution pair exists for each $\clZ_T$-measurable
$F$~\cite[Ch.~7]{yong1999stochastic}.  

\subsection{A word of caution}

Although the interpretation is not used in this paper, the
formula in~\eqref{eq:wrong_backward}
is an example of the so called backward-in-time It\^{o} integral.  The integral is denoted as $\int_t^T (\dots)\stackrel{\leftarrow}{\ud Z_s}$ with a backward
arrow~\cite[Sec. 4.2]{nualart1988stochastic} .  Using the backward-in-time integral,
one defines a solution of~\eqref{eq:SDE_arrow} that is
backward-adapted to the filtration,
i.e., solution at time $t$ depends upon randomness from time $t$ to
$T$.  The backward-adapted solution is different from the forward-adapted solution of
the BSDE.  However, such constructions are referred
to as ``reverse time'' or ``time reversed'' form of
SDEs~\cite{anderson1982reverse,nualart1988stochastic}.  

Historically,
in the filtering context, the backward Zakai equation is an example of
a SDE with a backward-adapted solution~\cite{pardoux1981non,benevs1983relation}.  More
recently, SDEs with 
backward-adapted solutions have become popular for generative AI
applications~\cite{song2020score}.

A reader is cautioned that, for us, ``reversing the arrow of time'' means a 
solution of a BSDE.  In particular, all the stochastic processes 
considered in this paper are forward-adapted.  For example, in
formula~\eqref{eq:right_backward}, the integral ``$\int_t^T V_s \ud
Z_s$'' is a standard It\^{o}-integral where the integrand
$V_s \in \clZ_s$ for $t\leq s \leq T$. 
The BSDE formalism is {\em necessary} to
generalize the Kalman's PoD and duality taught in an introductory control systems
course. 

\clearpage 
\section{Mathematical foundations of duality}
\label{sidebar:math_found_duality}
\subsection{Dual vector spaces}\label{ssec:dual-space}

Let $\clY$ be a Banach space equipped with norm $\|\cdot\|_\clY$. The {\em dual space}, denoted by $\clY^\dagger$, is the space of bounded linear functionals on $\clY$.  For $y^\dagger\in\clY^\dagger$, the notation 
\[
\langle y,y^\dagger\rangle := y^\dagger(y)
\]
is used to denote the evaluation at $y\in\clY$. The bilinear map
$\langle \cdot ,\cdot\rangle:\clY\times \clY^\dagger \to \Re$ is
called the {\em duality pairing}.  The dual space $\clY^\dagger$ is
also a Banach space~\cite[Thm.~5.3.1]{luenberger1997optimization}. 
An important special case is when the space $\clY$ is a Hilbert
space equipped with an inner product $\langle \cdot,\cdot\rangle$.
For a Hilbert space, $\clY^\dagger = \clY$ and the duality pairing is
the inner product.  In case there is more than one function
space, the duality pairing is indicated by using a subscript, e.g,
$\langle \cdot ,\cdot\rangle_\clY$.

\subsection{Adjoint}

Consider another Banach space $\clU$ and let $\clL:\clU\to\clY$ be a
bounded linear operator. The adjoint operator
$\clL^\dagger:\clY^\dagger\to \clU^\dagger$ is obtained from the
defining relation: 
\[
\langle u, \clL^\dagger y^\dagger \rangle_\clU = \langle \clL u,
y^\dagger\rangle_\clY,\quad \forall y^\dagger \in \clY^\dagger,\;u\in \clU
\]
A detailed explanation showing that $\clL^\dagger$ is well-defined,
linear and bounded operator appears
in~\cite[Ch.~6.5]{luenberger1997optimization}.

\subsection{Closed range theorem}

In matrix theory, it is an elementary fact that the range space of a
matrix is orthogonal to the null space of its
transpose~\cite{strang1993fundamental}. The following theorem provides
a generalization of this fact to the relationship between the range space $\Rsp(\clL)$ of
$\clL$ and the null-space $\Nsp(\clL^\dagger)$ its adjoint
$\clL^\dagger$. 

\medskip

\begin{theorem}[Theorem 6.6.1 in~\cite{luenberger1997optimization}]\label{thm:dual-adjoint}
	Let $\clU$ and $\clY$ be Banach spaces and $\clL:\clU\to\clY$ is a bounded linear operator. Then
	\[
	\Rsp(\clL)^\bot = \Nsp(\clL^\dagger)
	\]
\end{theorem}

\subsection{Relationship to the present paper}

In this paper, the classical duality between controllability and
observability of deterministic linear systems is reviewed to be a
special case of the abstract duality in linear algebra.  The review
sets the stage to discuss a novel extension of the duality for HMMs. 

For both cases -- linear system and the HMM -- duality entails definition
of three objects: (a) function space of admissible control inputs $\clU$; (b)
function space of ``dual'' state $\clY$; and (c) a linear operator
$\clL:\clU\to \clY$.
For deterministic linear systems, these definitions are classical and reviewed here.  For the HMM, these definitions are intended such that controllability and observability are dual properties based on an application of the closed range theorem.

The operator $\clL$ is the solution operator of a input-state system which is referred to as the {\em dual control system}.  A novel contribution of our work is to define the dual control system for an HMM.

\newpage
\clearpage
\section{Derivation of duality principle for the linear Gaussian model}\label{sidebar:derivation_KF}

\revblue{
\subsection{Duality principle}

Consider the state process $X$ defined by~\eqref{eq:linear-Gaussian-model-a} and the dual process $y$ defined by~\eqref{eq:mv-intro-b}, respectively.  Then
\[
\ud (y_t^\tp X_t) = u_t^\tp \ud Z_t + u_t^\tp \ud W_t + y_t^\tp \sigma \ud B_t,\qquad 0\leq t\leq T
\]
Integrating both sides from $0$ to $T$, 
\[
f^\tp X_T - \Big(\underbrace{y_0^\tp m_0 - \int_0^T u_t^\tp \ud Z_t}_{=S_T}\Big) = \big(y_0^\tp X_0 - y_0^\tp m_0\big) + \int_0^T u_t^\tp \ud W_t + y_t^\tp \sigma\ud B_t
\]
Each of the terms on the right-hand side are mutually independent and have zero mean. Therefore, upon squaring and taking expectation,
\[
\E\big(|f^\tp X_T - S_T|^2\big) = y_0^\tp \Sigma_0 y_0 + \int_0^T |u_t|_R^2 + |y_t|^2_Q \ud t
\]
Therefore, the m.s.e of the estimator~\eqref{eq:LG-estimator} is expressed as an optimal control objective~\eqref{eq:mv-intro-a}.
}

\subsection{Derivation of the Kalman-Bucy filter}
The optimal solution to a linear-quadratic (LQ) problem~\eqref{eq:mv-intro} is given in a linear feedback form~\cite[Theorem 3.1]{bensoussan2018estimation}: 
\[
u_t = -R^{-1} H^\tp \Sigma_ty_t
\]
where $\Sigma_t$ is the solution to the (forward-in-time) dynamic Riccati equation (DRE):
\begin{equation}\label{eq:Ricc-LG}
	\frac{\ud }{\ud t}\Sigma_t = A^\tp \Sigma_t + \Sigma_t A + Q -
        \Sigma_t H R^{-1} H^\tp \Sigma_t,\quad \Sigma_0\text{ given}
\end{equation}
It was noted by Kalman and Bucy that the DRE~\eqref{eq:Ricc-LG} for
the LQ problem is identical to the
DRE~\eqref{eq:KF_var} of the Kalman filter~\cite{kalman1961}.  
Let $\Phi(T,t)$ be the transition matrix from time $T$ to $t$ of the closed-loop system
\[
-\frac{\ud }{\ud t}\Phi(T,t) = (A-HR^{-1} H^\tp \Sigma_t )\Phi(T,t),\quad \Phi(T,T) = I
\]
Substituting the optimal control into~\eqref{eq:LG-estimator},
\begin{align*}
	f^\tp \hat{X}_T &= y_0^\tp m_0 + \int_0^T y_t^\tp \Sigma_t H
                          R^{-1} \ud Z_t\\
	&= f^\tp\Phi^\tp(T,0) m_0 + \int_0^T f^\tp \Phi^\tp(T,t)
          \Sigma_t H R^{-1} \ud Z_t
\end{align*}
Since $f$ is arbitrary,
\[
\hat{X}_T = \Phi^\tp(T,0)m_0 + \int_0^T \Phi^\tp(T,t)\Sigma_t H R^{-1} \ud Z_t
\]
Because $T$ is arbitrary, denote it as $t$
\[
\hat{X}_t = \Phi^\tp(t,0)m_0 + \int_0^t \Phi^\tp(t,s)\Sigma_s H R^{-1} \ud Z_s
\]
Differentiating both sides with respect to $t$ yields the equation of the Kalman-Bucy filter:
\begin{align}
	\ud \hat{X}_t &= (A-H R^{-1} H^\tp\Sigma_t )^\tp\Big( \Phi^\tp(t,0)m_0 + \int_0^t 
	\Phi^\tp(t,s)\Sigma_s H R^{-1} \ud Z_s\Big)\ud t + \Sigma_t H R^{-1} \ud 
	Z_t \nonumber\\
	&= A^\tp \hat{X}_t \ud t + \Sigma_t H R^{-1} \big(\ud Z_t - H^\tp \hat{X}_t\ud t\big) 
	\label{eq:KF-equation}
\end{align}

\clearpage
\newpage
\section{Simple example to illustrate the two types of duality}
\label{sidebar:simple_example}

The simple example is adapted from~\cite[Section
3.5]{kailath2000linear} to illustrate the main ideas.
Consider a linear estimation problem defined by the model:
\begin{equation*}\label{eq:basic-problem}
	Z = H^\tp X_0 + W
\end{equation*}
where $X_0\sim N(m_0,\Sigma_0)$, $W\sim N(0,R)$ are independent
Gaussian random variables of dimension $d$ and $m$, respectively.  It
is assumed that both $\Sigma_0$ and $R$ are strictly positive-definite.
The goal is to compute the conditional mean $\E(X_0\mid Z)$.

\subsection{Minimum variance construction}

Fix $f\in \Re^d$. The estimation objective is to compute $\E(f^\tp X_0
\mid Z)$. 
Since all random variables are Gaussian, it suffices to consider an estimator 
$S$ of the form
\begin{equation}\label{eq:simple-estimator}
	S= b - u^\tp Z 
\end{equation}
where $b \in \Re$ and $u \in \Re^m$ are both
deterministic~\cite[Corollary 1.10]{le2016brownian}. 

The minimum
variance optimization problem is as follows:
\begin{equation}\label{eq:basic-stochastic-problem}
	\min_{\substack{b\in\Re, u\in\Re^m}}\E\big(|f^\tp X_0 - S|^2\big)
\end{equation}
With the estimator~\eqref{eq:simple-estimator}, the optimization objective 
becomes
\[
\E\big(|f^\tp X_0 - S|^2\big) = |f + H u|^2_{\Sigma_0} + |u|_R^2 + 
\big((f + H u)^\tp m_0 - b\big)^2
\]
Set $y = f+H u$.  It then follows that
$b = y^\tp m_0$ is the optimal choice, and the minimum 
variance optimization problem becomes a quadratic programming problem: 
\begin{align*}
	\min_{u\in\Re^m}\quad & |y|^2_{\Sigma_0}  + |u|_R^2 \\
	\text{s.t. } \quad&y = f+H u
\end{align*}
Its solution is given by
\[
u = (H^\tp \Sigma_0 H +R)^{-1} H^\tp \Sigma_0 f
\]
and the corresponding optimal estimator is
\[
S= f^\tp \big(m_0 + \Sigma_0 H(H^\tp \Sigma_0 H +R)^{-1}(Z-H^\tp m_0)\big) 
\]
Since $f$ is arbitrary, 
\begin{equation}\label{eq:solution-simple-example}
	\E(X_0\mid Z) = m_0 + \Sigma_0 H(H^\tp \Sigma_0 H +R)^{-1}(Z-H^\tp m_0)
\end{equation}

\subsection{Minimum energy/maximum likelihood construction}

While the previous problem considers the (minimum variance) property of the conditional expectation, the minimum energy problem begins with the Bayes' formula for conditional density:
\[
\rho_{X_0 \mid Z}(x\mid z) = \frac{\rho_{X_0,Z}(x,z)}{\rho_Z(z)},\quad x\in\Re^d, z\in \Re^m
\]
where $\rho_{X_0,Z}$ denotes the joint probability density function, $\rho_Z$ is the marginal and $\rho_{X_0\mid Z}$ denotes the conditional density.
The objective is to compute the maximum-likelihood estimate of $X_0$.
Since the event $[X_0=x,Z=z]$ is the same as $[X_0=x,W=z-Hx]$, we have
\[
-2\log\big(\rho_{X\mid Z}(x\mid z)\big) = |x-m_0|^2_{\Sigma_0^{-1}} +
|z-H^\tp x|^2_{R^{-1}} + c(z)
\]
where the constant $c(z)$ only depends on $z$. Therefore, the maximum 
likelihood problem is given by 
\begin{equation}\label{eq:basic-deterministic-problem}
	\min_{x\in \Re^d}\;  |x-m_0|^2_{\Sigma_0^{-1}} +
|z-H^\tp x|^2_{R^{-1}} 
\end{equation}
Its optimal solution is obtained as
\[
x = (\Sigma_0^{-1} + H R^{-1} H^\tp )^{-1}(\Sigma_0^{-1} m_0 + H R^{-1} z)
\]
By an application of the matrix inversion
lemma~\cite[Appdx. A.1]{kailath2000linear}, this formula is identical
to~\eqref{eq:solution-simple-example} with $Z=z$.

\medskip

	In~\cite{kailath2000linear},~\eqref{eq:basic-stochastic-problem}
        is referred to as the stochastic
        problem and~\eqref{eq:basic-deterministic-problem}
	 as the deterministic problem.  
	It is noted 
        that the objectives~\eqref{eq:basic-stochastic-problem}
        and~\eqref{eq:basic-deterministic-problem} are not directly
        related to each other even though they share the same
        solution~\cite[p.~100]{kailath2000linear}.

\subsection{Relationship to the present paper}

In this paper, the two types of dualities for the linear Gaussian
model~\eqref{eq:linear-Gaussian-model} are reviewed as generalizations
of the basic constructions described in this sidebar.  The review sets
the stage for describing their extensions for HMMs.  The
minimum variance extension is new based on original work of the
authors.

\newpage
\clearpage

\revblue{
\section{Minimum energy duality as a deterministic optimal control problem}

Concerning the dual optimal control problem~\eqref{eq:mee-intro},
Bensoussan writes in~\cite[p.~180]{bensoussan2018estimation}: ``{\em The notation is reminiscent of the probabilistic origin. The function $z$ is a given $L^2\big([0,T];\Re^m\big)$ function. It is reminiscent of the observation process, in fact rather the
			derivative of the observation process (which, as we know, does not exist). Similarly,
			$v$ is reminiscent of the noise that perturbs the system (again its derivative),
			and $x_0$ is the value of the initial condition, which we do not know. The cost
			functional~\eqref{eq:mee-intro} contains
                        weights related to the covariance matrices
                        that were part of the initial probabilistic
                        model.
		}''

\subsection{Deterministic least squares}
  The minimum energy optimal control objective~\eqref{eq:mee-intro-a} is ubiquitous in the estimation literature focussed on the deterministic settings of the problem.  It is often introduced as a starting point, without any reference to the stochastic filtering model.  For example, Willems in his paper~\cite{willems2004deterministic}, refers to~\eqref{eq:mee-intro} as ``deterministic least squares filtering''. The main result (Theorem~3) of~\cite{willems2004deterministic} is the derivation of~\eqref{eq:LG-smoother-forward} which he refers to as the Kalman filter; \cite[Sec.~8.3]{sontag2013mathematical} refers to this derivation as the ``(deterministic) Kalman filter''.  Willems writes: ``{\em The purpose of this paper is to show that the Kalman filter admits a perfectly satisfactory deterministic least squares formulation. I will argue that this approach is eminently reasonable, and circumvents the sticky modeling assumptions that are unavoidable in the stochastic approach}''. According to Willems, the ``sticky modeling assumptions'' relate to the modeling of dynamic processes as stochastic processes, specifically the use of the SDE~\eqref{eq:linear-Gaussian-model-b} to model the output process. In his view, instead of justifying~\eqref{eq:mee-intro-a} from the perspective of log-likelihood, it is more appropriate to interpret~\eqref{eq:mee-intro} as the `natural' least-square objective for estimation, which is a classical topic with a well-established history~\cite{swerling1971modern}.  For the same reason, Fleming refers to the nonlinear extension of~\eqref{eq:mee-intro} as ``deterministic nonlinear filtering~\cite{fleming1997deterministic}''.

 Although we do not agree with some of the modeling critiques presented in Willems' paper, we do concur with his and others' interpretation of~\eqref{eq:mee-intro} as a deterministic least-squares problem.  This interpretation is in fact discussed at length in the chapter on duality in~\cite[Ch.~15]{kailath2000linear} where~\eqref{eq:mee-intro} is referred to as an `equivalent problem' -- called as such because it yields the same solution as the stochastic problem. Kailath et.~al.~explicitly note that~\eqref{eq:mee-intro} is {\em not} a dual to the stochastic filtering problem~\cite[Sec.~15.3]{kailath2000linear}.

\subsection{Relationship between deterministic and stochastic problems}

An important research theme is to relate the solution (at the terminal time $T$) obtained from solving the deterministic (optimal control) and stochastic (filtering) problems. For the linear Gaussian model, both yield the equation of the Kalman filter (see formula~\eqref{eq:LG-smoother-forward}). In nonlinear settings, these questions were first discussed in Mortensen's paper~\cite[Sec.~5]{mortensen1968} and results have been obtained in the small-noise limit~\cite{hijab1980minimum,james1988nonlinear,fleming1997deterministic,krener2004convergence}.  These results are helpful to provide a justification for the extended Kalman filter as well as develop stability analysis for the same. 
}

\newpage
\clearpage

\section{Examples of HMM}

\subsection{Finite state-space HMM}

Let $\bS = \{1,2,\hdots,d\}$. The space of functions and measures are both identified with $\Re^d$: a real-valued function $f$ (resp., a measure $\mu$) is identified with a column vector in $\Re^d$ where the $x^{\text{th}}$ element of the vector represents $f(x)$ (resp., $\mu(x)$), for $x\in\bS$, and $\mu(f) = \mu^\tp f$.  In this manner, the
observation function $h:\bS\to\Re^m$ is also identified with a matrix
$H\in\Re^{d\times m}$.  $\clP(\bS)$ is the probability simplex in
$\Re^d$.   The generator is a transition rate matrix, denoted as
$A$, whose $(x,y)$ entry (for $x\neq y$) gives the non-negative rate of transition from 
$x\mapsto y$, for states $x,y\in\bS$.  The diagonal
entry $(x,x)$ is chosen such that the sum of the elements in the
$x$-th row is zero.   The finite state-space HMM is denoted as $(A,H)$.

\begin{figure}
  \centering
  \includegraphics[width=0.9\textwidth]{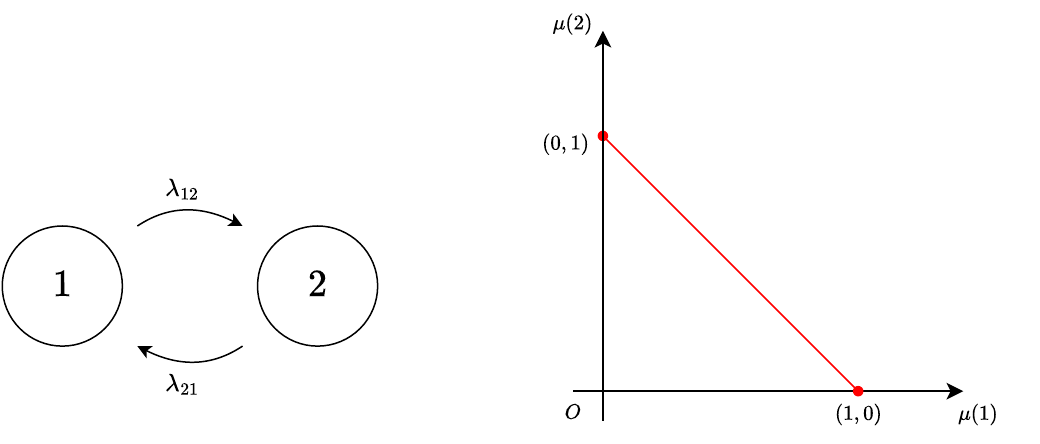}
  \caption{(a) State transition diagram for a 2-state Markov
    chain. (b) Probability simplex in the measure space $\Re^2$.} \label{fig:two_state}
\end{figure}

The simplest case is with $d=2$ and $m=1$.  Figure~\ref{fig:two_state} depicts the state-space, the transitions between the
two states, and the probability simplex.  The rate matrix and the
observation function are as follows:
\[
A = \begin{bmatrix} -\lambda_{12} & \lambda_{12} \\ \lambda_{21} &
  -\lambda_{21} \end{bmatrix},\qquad H =  \begin{bmatrix} h(1) \\ h(2) \end{bmatrix}
\] 

\subsection{It\^{o} diffusions}

An It\^o diffusion modeled using a
stochastic differential equation (SDE):
\begin{align}\label{eq:dyn_sde}
  \ud X_t& = a (X_t) \ud t + \sigma(X_t) \ud B_t,\quad X_0\sim\mu\\
  \ud Z_t& = h(X_t) \ud t + \ud W_t
\end{align}
where 
the drifts $a,h\in C^1(\Re^d;
\Re^d)$ and the diffusion coefficient $\sigma\in C^2(\Re^d; \Re^{d\times n})$ satisfy
appropriate technical conditions such that a strong solution exists
for $[0,T]$, and $B$ is a standard B.M. and $W$ is a B.M with covariance matrix $R\succ 0$. It is assumed that $X_0,B$ and $W$ are mutually independent.  With $\sigma=0$, the model for $X$ is an ODE.  

The infinitesimal generator $\clA$ acts on $C^2(\Re^d;\Re)$ functions in its domain according to~\cite[Thm. 7.3.3]{oksendal2003stochastic}
\begin{equation}\label{eq:generator_Ito_diffusion}
(\clA f)(x):= a^\tp(x) \nabla f(x) + \half  \tr\big(\sigma\sigma^\tp(x)(D^2f)(x)\big),\quad x\in\Re^d
\end{equation}
where $\nabla f$ is the gradient vector and $D^2 f$ is
the Hessian matrix.

\subsection{Linear Gaussian model}
The linear Gaussian model~\eqref{eq:linear-Gaussian-model} is a special case of an It\^o diffusion
  where the drift terms are linear $a(x)=A^\tp x$ and $h(x)=H^\tp x$,
  the diffusion coefficient $\sigma(x)=\sigma$ is a
  constant matrix, and the prior $\mu$ is a Gaussian density. 

A real-valued linear function is expressed as \[f(x) = 
\tilde{f}^\tp x,\quad x\in\Re^d\] 
where $\tilde{f}\in \Re^d$. Then it is an easy calculation to see that $\clA f$ is also a linear function given by
\begin{equation}\label{eq:gen-LG}
\big(\clA f\big)(x) = (A\tilde{f})^\tp x,\quad x\in\Re^d
\end{equation}

\newpage
\clearpage

\section{Duality between the HMM and BSDE control system}

\subsection{Background on nonlinear filtering}

A standard approach is based upon the Girsanov change of measure.  Let
$(\Omega,\clF_T,\sP)$be the probability triple and suppose $h\in C_b(\bS)$.  
Define a new probability measure $\tsP$ as follows:
\[
\frac{\ud \tsP}{\ud \sP} = \exp\Big(-\int_0^T 
h^\tp(X_t) \ud W_t - \half \int_0^T |h(X_t)|^2\ud t\Big) =: D_T^{-1}
\]
Then it is known that the probability law for $X$ is unchanged but
$Z$ is a $\tsP$-B.M.~that is independent of
$X$~\cite[Ch.~5]{xiong2008introduction}.  The
expectation with respect to $\tsP$ is denoted by $\tE(\cdot)$.  
The \emph{un-normalized filter} is a 
measure-valued process $\sigma = \{\sigma_t:0\le t \le
T\}$ defined by
\[
\sigma_t(f) := \tE\big(D_tf(X_t)|\clZ_t\big),\quad f\in 
C_b(\bS)
\]
Because $Z$ is a $\tsP$-B.M., the equation of unnormalized filter is
easily obtained and is in fact the celebrated Duncan-Mortensen-Zakai (DMZ) equation of
nonlinear filtering~\cite[Thm.~5.5]{xiong2008introduction}:
	\begin{equation}\label{eq:Zakai}
		\ud \sigma_t(f) = \sigma_t(\clA f) \ud t +
                \sigma_t(hf)^\tp \ud Z_t, \quad \sigma_0=\mu
	\end{equation}
The solution is denoted as $\sigma^\mu(f) := \{\sigma_t^\mu(f):0\leq t\leq T\}$ to stress the
dependence upon initial prior $\mu\in\clP(\bS)$.

\subsection{Function spaces}

The function spaces are as follows:
\begin{itemize}
\item $\clY = C_b(\bS)$ and its dual space $\clY^\dagger = \clM(\bS)$,
  the space of regular, bounded and finitely additive signed measures
  (rba measures).  The duality pairing
\[
\langle \mu, f \rangle_\clY := \mu (f) = \int_{\bS} f(x) \ud
\mu(x),\quad f\in  C_b(\bS),\;\mu\in \clM(\bS)
\]
\item The Hilbert space for $\Re^m$-valued stochastic processes 
\[
\clU = L^2_\clZ([0,T];\Re^m) := \Big\{U:\Omega\times[0,T]\to \Re^m: U\text{ is $\clZ$-adapted}\;\text{and}\; \tE\Big(\int_0^T|U_t|^2\ud t\Big)<\infty\Big\}
\]
equipped with the inner product
\[
\langle U, V\rangle_\clU := \tE\Big(\int_0^T U_t^\tp V_t \ud
t\Big),\quad U,V\in\clU
\]
\end{itemize}
The two function space $\clY$ and $\clU$ are counterparts of $\Re^d$
and $L^2([0,T];\Re^m)$ that are useful in the study of deterministic
linear systems~\eqref{eq:LTI-obs} and~\eqref{eq:LTI-ctrl}.  Note that for $\bS=\{1,2,\hdots,d\}$, the function space
$\clY=\Re^d$.  In all cases, $L^2([0,T];\Re^m)$ is a subspace of
$\clU$ comprising of deterministic signals. 

\subsection{Linear operator and its adjoint}

Consider the BSDE control system~\eqref{eq:dual-bsde} with control
input $U\in\clU$.  Its solution operator is used to define a bounded linear operator $\clL:
\clU\times \Re \to \clY$ as follows:
\begin{equation*}\label{eq:ctrl-operator}
	\clL(U,c) = Y_0,\quad U\in\clU,\;c\in\Re
\end{equation*}
where $Y_0\in \clY$ is the solution at time 0 (Note that $Y_0$ is
deterministic). This operator is the counterpart of the linear
operator, also denoted as $\clL$, previously defined for the
input-state 
system~\eqref{eq:LTI-ctrl}.

In~\cite[Thm.~2]{duality_jrnl_paper_I}, it is shown that the adjoint of $\clL$ is
the linear operator $\clL^\dagger: \clM(\bS)\to L^2_\clZ\big([0,T];\Re^m\big)\times \Re$ as follows:
\begin{equation}\label{eq:observability-operator}
	\clL^\dagger \mu = \big(\{\sigma_t^\mu(h):0\le t\le T\},
        \mu(\ones)\big) =: (\sigma^\mu(h), \mu(\ones))
\end{equation}
The proof entails showing the following duality relationship:
\[
\langle \mu,\clL(U,c) \rangle_\clY = \big\langle \sigma^\mu(h) , U \rangle_{\clU} + c
\mu(1) = \big\langle \clL^\dagger \mu, (U,c)\big\rangle_{\clU\times
  \Re},\quad \forall \; U\in\clU,\;\mu \in \clM(\bS)
\]

\subsection{Controllability, observability, duality}

Based on $\clL:
\clU\times \Re \to \clY$, controllability is defined in the same way as in
the linear systems theory. 

\begin{definition}[Defn.~3 in~\cite{duality_jrnl_paper_I}]
	For the BSDE~\eqref{eq:dual-bsde}, the \emph{controllable subspace} $\clC_T := \Rsp(\clL)$. 
Explicitly,
\begin{equation}\label{eq:ctrl-subspace}\clC_T = \big\{y_0 \in \clY: \exists\, c\in \Re\text{
  and } U \in \clU\; \text{s.t.}\; Y_0 = y_0\;\text{and}\; Y_T = c\ones\big\}
\end{equation}
	The BSDE~\eqref{eq:dual-bsde} is said to be \emph{controllable} if $\clC_T$ is dense in $\clY$.
\end{definition}

The following definition is adopted for observability of the HMM
(see~\cite[Thm.~1]{duality_jrnl_paper_I} where this definition is
shown to be equivalent to van Handel's
definition of stochastic observability~\cite[Defn.~3]{van2009observability}):

\begin{definition}
	The HMM $(\clA,h)$ is \emph{observable} if:
	\[
	\sigma_t^\mu(h) = \sigma_t^\nu(h), \quad 0\leq t\leq T \quad \implies \quad \mu = \nu,\quad
        \forall \;\;\mu,\nu \in \clP(\bS)
	\]
\end{definition}

\begin{figure}
	\centering
	\includegraphics[width=0.9\textwidth]{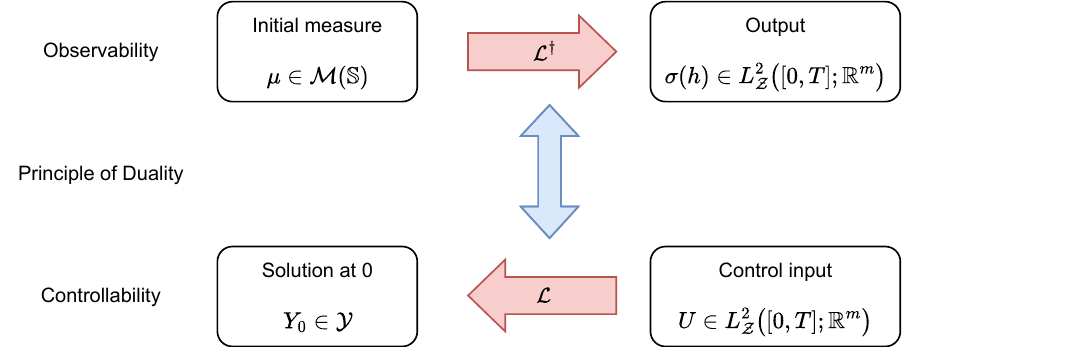}
			\caption{Duality for HMM: Controllability and observability are properties of a linear operator and its adjoint.}
		\label{fig:duality_NL}
\end{figure}

From closed range theorem, controllability and observability are now
seen to be dual properties.  Specifically, 
the dual control system is controllable iff HMM is observable.  This
is the promised generalization of the elementary duality between
observability of the state-output system~\eqref{eq:LTI-obs} and the
controllability of the input-state system~\eqref{eq:LTI-ctrl}. 
Figure~\ref{fig:duality_NL} depicts the duality relationship (compare with Fig.~\ref{fig:duality_LTI}). 

\revblue{
  \subsection{Proof sketch for duality principle~\eqref{eq:duality-principle} for HMM $(\clA,h)$}
The idea of the proof is the same as proof of the duality principle for the linear Gaussian model (see ``Derivation of duality principle for the linear Gaussian model''). 
  
All the stochastic processes, the state and observation process $(X,Z)$ and the dual process $(Y,V)$, are forward adapted.  Apply the It\^o-Wentzell formula to $Y_t(X_t)$ (see~\cite[Thm.~1.17]{rozovskiui2018stochastic} for the formula in the Euclidean case and~\cite{kunita1981some} for
the general case) to obtain,
\begin{align*}
	\ud Y_t(X_t) 
	&= - U_t^\tp \ud Z_t + \big(U_t + V_t(X_t)\big)^\tp \ud W_t + \ud N_t
\end{align*}
where $N=\{N_t:0\leq t\leq T\}$ is a certain martingale associated with the Markov process (see~\cite[Remark~1]{duality_jrnl_paper_II}). Define the error $\varepsilon:=\{\varepsilon_t:0\leq t\leq T\}$ by
\[
\varepsilon_t := Y_t(X_t) - (\mu(Y_0) - \int_0^t U_s^\tp \ud Z_s),\quad 0\le t \le T
\]
whereby
\begin{align}\label{eq:varepsilon_DP}
\ud \varepsilon_t &= \big(U_t + V_t(X_t)\big)^\tp \ud W_t + \ud N_t
\end{align}
with $\varepsilon_0 = Y_0(X_0)-\mu(Y_0)$. 
Using the It\^o-rule,
\[
\ud \varepsilon_t^2 = 2\varepsilon_t \ud \varepsilon_t + |U_t+V_t(X_t)|_R^2 \ud t + \ud \langle N\rangle_t
\]
where $\langle N\rangle=\{\langle N\rangle_t:0\leq t\leq T\}$ is the quadratic variation of the process $N$. The formula for its differential is given by (see~\cite[Eq.~(6)]{duality_jrnl_paper_II}),
\[
\ud \langle N\rangle_t = (\Gamma Y_t)(X_t)\ud t
\]
where $\Gamma$ is the carr\'e du champ. 
Upon using this formula, the duality principle~\eqref{eq:duality-principle} is obtained by integrating~\eqref{eq:varepsilon_DP} and taking expectation because $\varepsilon$ is a martingale and $\E\big(|F(X_T)-S_T|^2\big) = \E\big(\varepsilon_T^2\big)$.  
}

\newpage
\clearpage 

\section{Carr\'e du champ operator}
\label{sidebar:carreduchamp}

For a Markov process with generator $\clA$, the \emph{carr\'e du champ} operator
        $\Gamma$ is defined as follows:
\[
\Gamma(f, g)(x) := (\clA fg)(x) - f(x)(\clA g)(x) - g(x)(\clA
f)(x),\quad x\in\bS
\]
where $f,g$ are
elements of a suitable space of test functions (denoted as $\clD$) such that 
each of the three terms on the righthand-side are meaningfully
defined~\cite[Defn.~3.1.1]{bakry2013analysis}.  $(\Gamma f)(x):= \Gamma(f,f)(x)$.   

For the models introduced in the sidebar ``Examples of HMM'', the formulae for the carr\'e du champ are as follows:
\begin{itemize}
\item Finite-state HMM $(A,H)$
\begin{equation}\label{eq:Gamma-finite}
	(\Gamma f)(x) = \sum_{y \in \mathbb{S}} A(x,y) (f(x) - f(y))^2,\quad x\in\bS
\end{equation}
\item For the 2-state example,
\[
(\Gamma f) =\begin{bmatrix} 1\\1 \end{bmatrix}  (\lambda_{12} + \lambda_{21}) (f(1)-f(2))^2
\]
\item For the It\^{o} diffusion~\eqref{eq:dyn_sde}, for a function $f\in C^1(\Re^d;\Re)$,
\begin{equation}\label{eq:Gamma-Euclidean}
	(\Gamma f) (x) = \big|\sigma^\tp(x) \nabla f(x) \big|^2,\quad x\in\Re^d
\end{equation}
\item For the linear Gaussian model~\eqref{eq:linear-Gaussian-model}, a real-valued linear function is expressed as $f(x) = \tilde{f}^\tp x$ with $\tilde{f}\in \Re^d$ and $\Gamma f$ is a constant function given by
\begin{equation}\label{eq:Gamma-LG}
\big(\Gamma f\big)(x) = \tilde{f}^\tp \big(\sigma \sigma^\tp\big) \tilde{f},\quad x\in\Re^d 
\end{equation}
\end{itemize}

\newpage
\clearpage

\section{Historical survey of minimum energy duality for nonlinear estimation}
\label{sidebar:history}

For the problems of nonlinear filtering and smoothing, solution approaches in literature based on duality include the following:
\begin{itemize}
	\item Mortensen's maximum likelihood nonlinear filter~\cite{mortensen1968}.
	\item Minimum energy estimator (MEE) such as the full
          information estimator (FIE) and the moving horizon estimator
          (MHE)~\cite[Chapter 4]{rawlings2017model}.
	\item Fleming-Mitter duality relating the Zakai equation of
          nonlinear filtering to the 
          Hamilton-Jacobi-Bellman (HJB) equation of optimal control~\cite{fleming1982optimal}.
	\item Mitter-Newton's variational formulation of nonlinear smoothing~\cite{mitter2003}.
\end{itemize}

A common theme connecting all of these prior works is that they are
all variations and generalizations of the minimum energy
duality~\eqref{eq:mee-intro} given therein for the linear Gaussian
model~\eqref{eq:linear-Gaussian-model}. While there are 
differences in specification of the optimal control objective, the
constraint in all these cases is a modified copy of the signal
model. Additional details on each of the four approaches
appears in the following subsections.

\subsection{Mortensen's maximum likelihood nonlinear filter}

In a pioneering paper, Mortensen~\cite{mortensen1968} considered the
maximum likelihood estimation (smoothing) problem for the following
model (compare with~\eqref{eq:linear-Gaussian-model}):
\begin{subequations}\label{eq:MLE-model}
\begin{align}
	\ud X_t &= a(X_t)\ud t + \sigma(X_t) \ud B_t,\quad X_0 \sim \mu \\
	\ud Z_t &= h(X_t)\ud t + \ud W_t,\quad Z_0 = 0
\end{align}
\end{subequations}
where modeling for the process noise $B$, the observation noise $W$, and the prior $\mu$ is as before, described for the linear Gaussian model~\eqref{eq:linear-Gaussian-model}.  Mortensen's objective was to
compute the maximum likelihood (ML) trajectory $x = \{x_t\in\Re^d:0\le t
\le T\}$ which is defined as the trajectory that maximizes the conditional density
\[
\rho_{X\mid Z}(x\mid z)
\]
given the output $z = \{z_t\in\Re^m: 0\le t \le T\}$.  

\begin{itemize}
\item Maximum likelihood estimation (MLE) problem (compare with~\eqref{eq:mee-intro})
\end{itemize}
\begin{subequations}\label{eq:MLE-problem}
	\begin{align}
		\mathop{\text{Minimize:}}_{x_0\in\Re^d, v\in L^2([0,T];\Re^p)}\bsJ_T(u,x_0;\dot{z}) &= |x_0-m_0|^2_{\Sigma_0^{-1}} + \int_0^T |v_t|^2 + |\dot{z}_t-h(x_t)|_{R^{-1}}^2 \ud t \label{eq:MLE-problem-a}\\
		\text{Subject to:}\qquad\qquad \;\;	\frac{\ud x_t}{\ud t} &= a(x_t) + \sigma (x_t) v_t \label{eq:MLE-problem-b}
	\end{align}
\end{subequations}
Note that the dimension of the control input $v$ is the same as the
process noise $B$ and that the objective function reduces to~\eqref{eq:mee-intro} for the linear
Gaussian model~\eqref{eq:linear-Gaussian-model}.  Similarity to~\eqref{eq:mee-intro}  comes from the assumed additive
Gaussian and independent nature of the noises $B$ and $W$, and the initial condition $X_0$ (see~\cite[Sec.~3.4.4]{JinPhDthesis} for a derivation). In Mortensen's paper, an algorithm to solve the MLE problem~\eqref{eq:MLE-problem} is proposed based on an application of the maximum principle to obtain the Hamilton's equations~\cite[Eqs.~(14)-(15)]{mortensen1968}.  The algorithm requires a forward and backward recursion to obtain the maximum likelihood trajectory.

Since Mortensen's early work, closely related optimization-type problem formulation have appeared for a plethora of filtering and smoothing problems. In different communities, these are referred by different names, e.g., maximum likelihood estimation (MLE), maximum a posteriori (MAP) estimation, minimum energy estimation (MEE), moving horizon estimation (MHE), and deterministic nonlinear filtering.

\subsection{Minimum energy estimator (MEE)}

While the Hamilton's equations can be used to obtain the optimal trajectory, it is also of interest to write a differential equation for the optimal filter -- nonlinear counterpart for the Kalman filter formula~\eqref{eq:LG-smoother-forward}.  The minimum energy estimator (MEE) is defined as the state (denoted as $\hat{x}_T$ at the terminal time $T$) that is an arg min of a certain value function (also at time $T$) for the optimal control problem~\eqref{eq:MLE-problem}.  The value function is defined by fixing the state $x_T=x$ and considering the minimum value of~\eqref{eq:MLE-problem-a} subject to the constraint~\eqref{eq:MLE-problem-b}.  Mortensen's goal was to write a (forward-in-time) differential equation for $\{\hat{x}_T:T\geq 0\}$.  Under the supervision of Art Krener, Hijab wrote an influential PhD thesis on this topic~\cite{hijab1980minimum} which was followed by a number of important works in this area~\cite{james1988nonlinear,fleming1997deterministic,krener2004convergence}.  For the linear Gaussian model, an accessible account of the MEE appears in the unpublished notes of~\cite{henderson2021tutorial}. 

\subsection{Minimum horizon estimator (MHE)}

Given the enormous success of model predictive control (MPC), related algorithms have been developed to solve the state estimation problems~\cite[Ch.~4]{rawlings2017model}.
In continuous-time setting, the optimal control problem is precisely the problem~\eqref{eq:MLE-problem}, typically augmented with additional constraints~\cite{goodwin2005}. Broadly, there are two classes of such algorithms: 
\begin{itemize}
	\item Full information estimator (FIE) where the entire history of observation is used.
	\item Moving horizon estimator (MHE) where the most recent window of observation is used.
\end{itemize}
We refer the reader to~\cite[Sec.~4.7]{rawlings2017model} where a
discussion on history of these approaches is provided. In this
section, the authors note that dual constructions are useful for
stability analysis. The authors describe certain results,
e.g.~\cite[Theorem 4.10]{rawlings2017model}, 
based on i-IOSS (incremental input/output to state stability) properties of the model.

\subsection{Fleming-Mitter-Newton duality}\label{ssec:MN-duality}
In~\cite{fleming1982optimal}, it is shown that the linear Zakai equation of
stochastic filtering can be transformed into the nonlinear Hamilton-Jacobi-Bellman
(HJB) equation of stochastic optimal control.
The particular transformation is an example of the log
transformation whereby the negative log of the posterior density is the value
function for a certain optimal control problem.
The precise form of the optimal control problem first appears in the 2003 paper of Mitter and Newton~\cite{mitter2003}. 
In this paper, the authors consider a control-modified version of the Markov process $X$ denoted by $\tilde{X} :=
\{\tilde{X}_t:0\le t\le T\}$. 
The control problem is to pick (a) the initial distribution
$\tilde\mu$ and (b) the state transition, such that
the law of $\tilde{X}$ equals the conditional law of $X$.  In essence, this is very similar to Mortensen's approach: whereas Mortensen aimed to find a single maximum likelihood trajectory, Mitter-Newton idea was to identify the entire stochastic process.

Abstractly, an optimization problem is formulated on the space of probability laws. Let $\sP_X$ denote the law for $X$, $\sQ$ denote the law for
$\tilde{X}$, and $\sP_{X\mid z}$ denote the law for $X$ given an observation path
$z=\{z_t:0\le t\le T\}$.  Assuming
$\sQ\ll\sP_X$, the objective function is the relative entropy (K-L
divergence) between $\sQ $ and $\sP_{X\mid z}$:
\begin{equation*}
	\min_{\sQ} \quad \E_{\sQ}\Big(\log \frac{\ud \sQ}{\ud \sP_X}\Big) - \E_{\sQ}\Big(\log\frac{\ud \sP_{X\mid z}}{\ud \sP_X}\Big).
\end{equation*}

\medskip

For the model~\eqref{eq:MLE-model}, this procedure yields the following
stochastic optimal control problem:
\begin{subequations}\label{eq:opt-cont-sde-hjb-intro}
	\begin{align}
		\mathop{\text{Min }}_{\tilde{\mu}, \; U}: \quad \sJ(\tilde{\mu},U\,;z) 
		& = \E\Big(\log \frac{\ud \tilde\mu}{\ud \mu}(\tilde{X}_0) + \int_0^T \ell(\tilde{X}_t,U_t\,;z_t)\ud t\Big)\\
		\text{Subj.} : \;\;\quad\qquad \ud \tilde{X}_t &= a(\tilde{X}_t)\ud t +
		\sigma(\tilde{X}_t)(U_t\ud t +
		\ud \tilde{B}_t), \quad
		\tilde{X}_0 \sim \tilde\mu
	\end{align}
\end{subequations}
where 
\[
l(x,u\,;z_t) := \half |u|^2 + \half |\dot{z}_t-h(x)|_{R^{-1}}^2
\]
The formula for cost function helps one see the minimum energy nature
of Mitter-Newton duality.   Because the derivative is only formally
defined, there is also a more rigorous
manner to express the cost function with only $z_t$ without the need
for $\dot{z}_t$~\cite[Rem.~2.2.5]{van2006filtering}.

The solution of the optimal control
problem~\eqref{eq:opt-cont-sde-hjb-intro} is given in the following
proposition which provides an explanation for the log transformation
originally noted in~\cite{fleming1982optimal}. 

\begin{proposition}
	Consider the optimal control problem~\eqref{eq:opt-cont-sde-hjb-intro}.
	For this problem, the HJB equation for the value function $\clV$ is
	as follows:
	\begin{align*}
		-\frac{\partial \clV_t}{\partial t}(x) &= \big(\clA (\clV_t+z_th)\big)(x) + \half h^2(x) -\half|\sigma^\tp(x)\nabla (\clV_t+z_th)(x)|^2,\quad0\leq t\leq T\\
		\clV_T(x) &= - z_Th(x),\quad x\in \Re^d
	\end{align*}
        where $\clA$ is the generator~\eqref{eq:generator_Ito_diffusion}. 
	The optimal control is of the state feedback form given by
        \[
        U_t^\opt = -\sigma^\tp(\tilde{X}_t) \nabla(\clV_t + z_th)(\tilde{X}_t),\quad 0\leq t\leq T
        \]
\end{proposition}

\medskip

\subsection{Log transform and time reversal}

Expressing $\clV_t(x) = -\log \big(q_t(x)e^{Z_th(x)}\big)$ it is readily verified
that $\{q_t:0\leq t\leq T\}$ solves the backward Zakai equation as follows:
\begin{align*}
-\ud q_t(x) &= {(\cal A} q_t)(x) \ud t + h(x)
	q_t(x) \overleftarrow{\ud Z}_t, \quad 0\leq t\leq T \\
	q_T (x) &= 1,\quad x\in\Re^d
\end{align*}
This reveals the log transform link between the HJB equation and the Zakai equation first described in~\cite{fleming1982optimal}.

The backward Zakai equation is also an example of a time reversal. Compared to a BSDE, the nature of time reversal in a backward Zakai is straightforward.  The crucial difference is that $\{q_t:0\leq t\leq T\}$ is a backward-adapted process (That is, its value at time $t$ depends upon randomness for times between $t$ to $T$).  In contrast, the solution of the BSDE is forward-adapted.  As explained in this paper, adaptedness of stochastic processes is an important property for the filtering problem.

\subsection{Additional reading}

A tutorial style review of the log transformation, its link to the
Zakai equation, specifically its path-wise robust representation,
formulations of the optimal control problem and its link to the
smoothing problem can be found in our review
paper~\cite{kim2020smoothing}. In addition to the Euclidean case,
explicit formulae are also described for the finite state-space.  
In~\cite[VI-B]{duality_jrnl_paper_II}, it is shown that, for the linear
Gaussian model~\eqref{eq:linear-Gaussian-model}, Mitter-Newton's problem~\eqref{eq:opt-cont-sde-hjb-intro} reduces to classical form of
minimum energy duality~\eqref{eq:mee-intro}.

Mitter-Newton duality was the basis for van Handel's PhD thesis on nonlinear filter stability~\cite{van2006filtering}.  Incidentally, the PhD thesis remains an isolated work, with van Handel himself moving on to an alternate non-duality approach (the so called intrinsic approach) in his later seminal papers on the topics of stochastic observability and filter stability~\cite{van2009observability,van2010nonlinear}.
Mitter-Newton duality has influenced algorithmic approaches for sampling and smoothing~\cite{sutter2016variational,reich2019data,chen2016relation,ruiz_kappen2017,taghvaei2023survey}.  More recent applications to topical problems and algorithms in machine learning are discussed in~\cite{raginsky2024variational}. 

\end{document}